\documentclass[11pt,a4paper]{amsart}

\usepackage{amsfonts, amssymb}
\usepackage{tikz}
\usepackage{tikz-cd}
\usepackage[margin=2.5cm]{geometry}
\usepackage{booktabs}
\usepackage{hyperref}
\usepackage{xfrac}
\usepackage{epsfig}
\usepackage{setspace}
\usepackage{enumerate}

\newcommand{\Aut}{\mathrm{Aut}}

\newcommand{\coker}{\mathrm{coker}}

\newcommand{\Eff}{\mathrm{Eff}}
\newcommand{\Ext}{\mathrm{Ext}}

\newcommand{\Hom}{\mathrm{Hom}}

\newcommand{\Ima}{\mathrm{Im} \,}

\newcommand{\MHS}{\mathrm{MHS}}

\newcommand{\NE}{\mathrm{NE}}
\newcommand{\Nef}{\mathrm{Nef}}

\newcommand{\Pic}{\mathrm{Pic}}

\newcommand{\rk}{\mathrm{rk} \,}

\newcommand{\cD}{\mathcal{D}}

\newcommand{\cO}{\mathcal{O}}

\newcommand{\cY}{\mathcal{Y}}

\newcommand{\CC}{\mathbb{C}}

\newcommand{\PP}{\mathbb{P}}
\newcommand{\QQ}{\mathbb{Q}}
\newcommand{\RR}{\mathbb{R}}

\newcommand{\ZZ}{\mathbb{Z}}
\newcommand{\Newt}{\mathrm{Newt}}

\newcommand{\ip}{\raise1pt\hbox{\large $\lrcorner$}}

\theoremstyle{plain}
\newtheorem{theorem}{Theorem}[section]
\newtheorem{proposition}[theorem]{Proposition}
\newtheorem{lemma}[theorem]{Lemma}

\newtheorem{corollary}[theorem]{Corollary}
\newtheorem{conjecture}[theorem]{Conjecture}
\newtheorem*{problem*}{Problem}
\newtheorem*{conjecture*}{Conjecture}
\newtheorem*{theorem*}{Theorem}

\theoremstyle{definition}
\newtheorem{definition}[theorem]{Definition}
\newtheorem{remark}[theorem]{Remark}

\title{A generic Torelli theorem for certain log Calabi--Yau threefolds}
\author{Wendelin Lutz}
\address{Leibniz University Hannover, Institut of Algebraic Geometry, Welfengarten 1,
	30167 Hannover, Germany}
\email{lutz@math.uni-hannover.de}
\begin{document}
	\begin{abstract}
		We prove a generic Torelli theorem for families $\pi \colon (\cY, \cD) \rightarrow B$ of smooth three-dimensional log Calabi--Yau pairs, under the assumptions that one fiber $(Y, D)$ arises as the blowup of a toric pair, and $H^3(Y,\ZZ)=0$. In the context of the Fano-LG correspondence, our result applies to very general fibers of families of Landau--Ginzburg models associated to a smooth Fano threefold with very ample anticanonical bundle.
	\end{abstract}
	\maketitle
	\section{Introduction}
	Throughout this paper, a \emph{log Calabi--Yau pair} $(Y, D)$ is a smooth projective variety $Y$ over $\CC$ with an anticanonical reduced, simple normal crossings divisor $D$. The quasiprojective variety $U=Y \setminus D$ is called a \emph{log Calabi--Yau variety}. 
	Log Calabi--Yau varieties share many properties of their compact cousins, but are often more amenable to explicit calculations: Many examples arise as iterated blowups of a toric pair $(\bar{Y}, \bar{D})$, so that combinatorial methods to study $(Y, D)$ are often available. This was used by Gross--Hacking--Keel to prove a global Torelli theorem for log Calabi--Yau \emph{surfaces} \cite{GrossHackingKeel}, and to prove versions of mirror symmetry for log Calabi--Yau surfaces \cite{GrossHackingKeelMS, HackingKeating}. \\
	Just as in the compact case, it is expected that some version of a Torelli theorem holds for log Calabi--Yau varieties.
	A log Calabi--Yau variety $U$ is quasi-projective, and therefore carries a mixed Hodge structure on its cohomology. A \emph{very} optimistic conjecture for a Torelli theorem might take the following form:
	\begin{conjecture}\label{conj:torelli1}
		Let $(\cY, \cD) \rightarrow B$ be a smooth family of $n$-dimensional log Calabi--Yau pairs. Let $b_1, b_2 \in B$ be closed points and suppose that parallel transport along some path $\gamma$ from $b_1$ to $b_2$ (in the analytic topology) induces an isomorphism 
		\[
		\gamma_* \colon H^n(\cY_{b_1} \setminus \cD_{b_1}, \ZZ) \cong H^n(\cY_{b_2} \setminus \cD_{b_2}, \ZZ)
		\]
		of mixed Hodge structures.
		Then there exist an isomorphism of pairs $(\cY_{b_1} , \cD_{b_1}) \cong (\cY_{b_2} , \cD_{b_2})$.
	\end{conjecture} 
	For our definition of families of pairs and parallel transport see Definition~\ref{def:paralleltransport}. 
	In dimension two, Conjecture~\ref{conj:torelli1} was proved by Gross--Hacking--Keel~\cite{GrossHackingKeel} if $D$ is singular, and by McMullen~\cite{McMullen} if $D$ is smooth.
	In dimension three and above, little is known. 
	Recently, Aguilar--Green--Griffiths~\cite{Griffiths} proved a Torelli theorem for $Y$ a generic cubic threefold, and $D$ a generic $K3$-surface.\\
	In this paper, we prove Conjecture~\ref{conj:torelli1} for a class of three-dimensional log Calabi--Yau pairs $(Y, D)$ with a \emph{toric model} (see Definition~\ref{def:toricmodel} and Lemma~\ref{lem:centersoftoricmodel}). Roughly speaking, this means that $Y$ is obtained as an iterated blowup $Y \rightarrow \dots \rightarrow \bar{Y}$ with $\bar{Y}$ a smooth toric variety, and $D$ is the strict transform of the toric boundary.
\begin{theorem}\label{thm:GenericTorelli}(Corollary~\ref{cor:main})
	Let $(\cY, \cD) \rightarrow B$ be a smooth family of log Calabi--Yau threefolds. Under the assumptions of Conjecture~\ref{conj:torelli1}, suppose further that 
	\begin{enumerate}
		\item $(\cY_{b_1} , \cD_{b_1})$ admits a toric model $(\bar{Y}, \bar{D})$, and every $1$-dimensional center in the toric model is a smooth rational curve. 
		\item The points $b_1, b_2 \in B$ are very general.
	\end{enumerate}
	Then there is an isomorphism of pairs $f \colon (\cY_{b_1} , \cD_{b_1}) \cong (\cY_{b_2} , \cD_{b_2})$ inducing $\gamma_*$.
\end{theorem}
Theorem~\ref{thm:GenericTorelli} includes a large class of interesting varieties, including families of log Calabi--Yau pairs which arise as mirrors to smooth Fano threefolds under the Fano/LG correspondence, see Section~\ref{sec:mirrorsymmetry}.\\
We make a few remarks about the assumptions in our theorem. 
It would perhaps be more natural to consider the class of all smooth log Calabi--Yau pairs $(Y, D)$ with \emph{maximal boundary}, by which we mean that $D$ admits a zero-dimensional stratum. If $(Y, D)$ admits a toric model, then it has maximal boundary. 
In dimension two, the converse also holds: a log Calabi--Yau surface pair has maximal boundary if and only if it admits a toric model (see \cite[Proposition 1.3]{GrossHackingKeelMS}).
This is no longer true in dimension three and above, see \cite{Kaloghiros} and references therein for several counterexamples. The existence of the toric model was key in the proof of the Torelli theorem of Gross--Hacking--Keel~\cite{GrossHackingKeel}, and is also a necessary ingredient in our proof.\\
We explain the condition that the only $1$-dimensional centers of the blowup are smooth rational curves: since $\bar{Y}$ is a smooth toric variety, it has no odd cohomology (\cite[Theorem~12.3.11]{CoxLittleSchenck}), and in particular, $H^3(\bar{Y})=0$. Denoting $Y=\cY_{b_1}$, the blowup formula for cohomology (see \cite[Theorem~7.31]{Voisin}), shows that $H^3(Y)=0$ if and only if every one-dimensional center in the toric model is a smooth rational curve. Equivalently, the weight $3$ part of the mixed Hodge structure $H^3(U)$ vanishes, so that $H^3(U)$ is a quotient of the mixed Hodge structure $H^2(D)$ of the (singular) boundary. This observation is crucial to our proof.\\
Just as for two-dimensional log Calabi--Yau pairs, we expect that the theorem remains true if the genericity assumption is dropped. However, in view of the weak Torelli theorem of \cite[Theorem 1.8]{GrossHackingKeel}, it is then likely no longer true that the isomorphism of pairs induces the isomorphism of mixed Hodge structures. We plan to address this in future work.
\subsection{Remarks on the statement of Global Torelli}
Let $Y$ and $Y'$ be deformation-equivalent, $n$-dimensional smooth projective Calabi--Yau manifolds. We say that $Y$ and $Y'$ are \emph{Hodge-equivalent} if there exists an isomorphism $H^n(Y, \ZZ) \cong H^n(Y', \ZZ)$ of polarized Hodge structures. 
The classical Torelli problem asks whether two Hodge-equivalent Calabi--Yau $n$-folds are isomorphic. This is known to be true in dimension two, and also in some three-dimensional examples: Voisin~\cite{Voisin2} proved a generic Torelli theorem for the quintic threefold, and Sheng--Xu~\cite{ShengXu} proved a global Torelli theorem for Calabi--Yau threefolds which arise as cyclic triple covers of $\PP^3$ branched along stable hyperplane arrangements. In general, the global Torelli theorem fails in dimension three: Szendr\H{o}i~\cite{Szendroi} constructed deformation-equivalent Calabi--Yau threefolds $Y$ and $Y'$ which are Hodge-equivalent and birational, but not isomorphic. Borisov--C\u{a}ld\u{a}raru--Perry~\cite{BorisovPerry} produced examples of Calabi--Yau threefolds which are not birational, but Hodge-equivalent modulo torsion. Ottem--Rennemo~\cite{Ottem} produced examples of Calabi--Yau threefolds which are not birational, but Hodge-equivalent. \\
In view of this, it makes sense to consider a weakened version of the Torelli problem which requires that the isomorphism of Hodge structures is induced by parallel transport in a smooth family $\cY \rightarrow B$ of Calabi--Yau manifolds, for some choice of path $\gamma$ in $B$ from $b_1$ to $b_2$, with $Y \cong \cY_{b_1}$ and $Y'=\cY_{b_2}$. To our knowledge, there are no counterexamples to this modified Torelli problem. Moreover, Verbitsky's Torelli theorem for Hyperk\"ahler manifolds~\cite{Verbitsky} relied on the assumption that the isomorphism of Hodge structures is induced by parallel transport. 
\subsection{Outline}
In Section 2, we review material on log Calabi--Yau pairs and the Torelli theorem for Looijenga pairs (log Calabi--Yau surfaces with maximal boundary).\\
In Section 3, we prove a Torelli theorem for type III surfaces, by which we mean a simple normal crossing surface $D$ with $\omega_D \cong \cO_D$ and dual complex $\Sigma_D$ homeomorphic to $S^2$ (Theorem~\ref{thm:typeIIITorelli}). Such surfaces appear as central fibers of Kulikov degenerations of $K3$-surfaces, so that this result may be of independent interest.
Since each component of $D$ is a Looijenga pair, we deduce this result from the Torelli theorem for Looijenga pairs and glueing. \\
Section 4 is devoted to the proof of Theorem~\ref{thm:GenericTorelli}. We give a brief synopsis: keeping the notation of the theorem, the toric model $(\cY_{b_1}, \cD_{b_1}) \rightarrow (\bar{Y}, \bar{D})$ is a sequence of contractions of $K$-negative extremal rays of the cone of curves (Mori contractions). We first show that $(\cY_{b_2}, \cD_{b_2})$ can be blown down to $(\bar{Y}, \bar{D})$ via a sequence of contractions of the ``same type". To achieve this, we use that $b_1, b_2 \in B$ are very general, so that parallel transport along the path $\gamma$ preserves the nef cone (see Lemma~\ref{lem:coneparalleltransport}). Dually, this allows us to match up $K$-negative extremal rays of the cone of curves, and we construct the blowdown $(Y', D') \rightarrow (\bar{Y}, \bar{D})$ using Mori's classification of $K$-negative extremal rays. \\
Having reduced to the case where $(\cY_{b_1}, \cD_{b_1})$ and $(\cY_{b_2}, \cD_{b_2})$ have a common toric model,
we show that the extension class of the MHS $H^3(U)$ encodes data on the centers of blowups of the toric model. Since we have an isomorphism of mixed Hodge structures, the blowup centers are the same, so that $(\cY_{b_1}, \cD_{b_1})$ and $(\cY_{b_2}, \cD_{b_2})$ are isomorphic.
\subsection{Acknowledgements}
I would like to thank Alessio Corti, Paul Hacking, Matt Kerr, Yu--Shen Lin, and Yan Zhou for helpful discussions.
\section{Preliminaries}
\subsection{Log Calabi--Yau pairs}
A \emph{log Calabi--Yau pair} (or log CY pair) $(Y, D)$ is a smooth projective variety $Y$ over $\CC$ with an anticanonical reduced, simple normal crossings divisor $D$. We recall that $D=~\sum_i D_i$ has simple normal crossings (snc) if each irreducible component $D_i$ is smooth, and the components intersect everywhere transversely. If $\dim(Y)=n$, we also call $(Y, D)$ a log CYn pair. 
A \emph{toric pair} $(\bar{Y}, \bar{D})$ is a (smooth) toric variety $\bar{Y}$ with $\bar{D}=\bar{Y} \setminus (\CC^\times)^n$ its toric boundary. 
\begin{definition}
	Let $(Y, D)$ be a log Calabi--Yau pair.
	\begin{itemize}
		\item 
		The \emph{dual complex} $\Sigma_D$ of $D$ (also denoted $\Sigma$ if $D$ is clear from context) is a CW complex whose vertices correspond to the irreducible components of $D$, and for every stratum $F \subset \bigcap_{i \in I} D_i$, we attach a $(|I|-1)$-dimensional cell. In the terminology of \cite{Hatcher}, $\Sigma_D$ is a regular $\Delta$-complex.
		\item We say that $D$ has \emph{maximal boundary} if $D$ has a zero-dimensional stratum. Equivalently, $\Sigma_D$ has maximal possible dimension $\dim(Y)-1$.
		\item Let $D_v$ be a codimension $i$-stratum of $D$, corresponding to a $i$-cell in $\Sigma_D$. The \emph{interior} of $D_v$, denoted $D_v^\circ$, is $D_v$ with all lower dimensional strata removed. 
	\end{itemize}
\end{definition}
The log CY pairs we consider arise as iterated blowups of a toric pair. We make the following definitions. 
\begin{definition}\label{def:toricmodel}
	Let $(Y,D)$ be a log CY pair. Let $f \colon Y' \rightarrow Y$ be the blowup of $Y$ along a smooth center $Z$. 
	\begin{itemize}
		\item We say that $f$ is a \emph{toroidal} (or \emph{corner}) blowup if the exceptional divisor $E$ has discrepancy $a(E, Y, D)=-1$. By \cite[Theorem 4.16]{KollarSingularities}, $Z$ is a boundary stratum of $D$. Denoting $D'$ the \emph{total} transform of $D$, the pair $(Y', D')$ is again log CY.
		\item We say that $f$ is an \emph{interior blowup} if $Z$ has simple normal crossings with $D$ (\cite[3.24]{Kollarsnc}) and the exceptional divisor $E$ has discrepancy $a(E, Y, D)=0$. 
		Denoting $D'$ the strict transform of $D$, the pair $(Y', D')$ is again log CY.
		\item 	 A \emph{ toric model} is a toric pair $(\bar{Y}, \bar{D})$ together with a birational morphism $f \colon (Y, D) \rightarrow (\bar{Y}, \bar{D})$ that factors as a composite of \emph{interior} blowups.
	\end{itemize} 
\end{definition}
The effect of toroidal and interior blowups on the dual complex is well-understood: if 
\[
f \colon (Y, D) \rightarrow (Y', D')
\]
is a toroidal blowup, $\Sigma_{D'}$ is a subdivision of $\Sigma_D$. If $f$ is an interior blowup, $\Sigma_{D'}$ and $\Sigma_{D}$ are identified as CW-complexes. (\cite[§9]{KollarXuFernex}).
\begin{remark}	
	In dimension two, our definition of a toric model is the standard definition of a toric model used for example in  \cite{GrossHackingKeel}.
	In dimension three, a toric model sometimes denotes any crepant birational map $f \colon (Y, D) \dashrightarrow (\bar{Y}, \bar{D})$ (e.g \cite[p.1]{Ducat}, \cite[Definition 2.6]{MoragaToricmodel}). In particular $f$, might not be a morphism and also include toroidal blowups. 
\end{remark}
We give an explicit description of the centers of an interior blowup in low dimension:
\begin{lemma}\label{lem:centersoftoricmodel}
	Let $(Y, D)$ be a log Calabi--Yau pair and let $(Y', D') \rightarrow (Y, D)$ be an interior blowup with smooth center $Z$ and exceptional divisor $E$. Then $Z \subset D$ and
	\begin{itemize}
		\item If $\dim(Y)=2$, then $Z$ is a point on the smooth locus of $D$.
		\item If $\dim(Y)=3$, then $Z$ is either a point in the interior of a one-dimensional stratum of $D$, or a smooth curve contained in a single irreducible component of $D$, meeting one-dimensional strata transversely, and avoiding zero-dimensional strata.
	\end{itemize}
	\begin{proof}
		Recalling that $Y$ is smooth and $D$ is simple normal crossing we have from the definition of discrepancy that 
		\[
		a(E, Y, D)=k-1-n
		\] 
		where $k$ is the codimension of $Z$ in $Y$ and $n$ is the number of irreducible components of $D$ containing $Z$.  The condition that $a(E, Y, D)=0$ forces $Z \subset D$. 
		If $\dim(Y)=2$, we see that $Z$ is a point on the smooth locus of $D$, and if $\dim(Y)=3$, then $Z$ is either a point in the interior of a one-dimensional stratum, or a curve contained in a single component of $D$.
		In the latter case, $Z$ has simple normal crossings with $D$ by assumption, so meets one-dimensional strata of $D$ transversely, and avoids zero-dimensional strata of $D$.
	\end{proof}
\end{lemma}
We finally define deformations of pairs and parallel transport.
\begin{definition}\label{def:paralleltransport} Let $\pi \colon \cY \rightarrow B$ be a smooth family of projective varieties over a connected analytic or algebraic base $B$. Let $\cD$ be a relative anticanonical divisor with simple normal crossings on $\cY$ (i.e. $\cD$ restricts to an anticanonical divisor on each fibre of $\pi$). 
	\begin{itemize}
		\item 	We say that $(\cY, \cD)$ is a \emph{family of log Calabi--Yau pairs} if the dual complex $\Sigma_{D_b}$ of fibers is locally constant (in the analytic topology). So $\Sigma_{D_b} \cong \Sigma_{D_b'}$ as CW complexes locally on $B$, 
		but the corresponding strata of $\cD$ might deform nontrivially.
		\item Let $(Y_1, D_1)$ and $(Y_2, D_2)$ be log Calabi--Yau pairs. If there are $b_1, b_2 \in B$ with $(\cY_{b_i}, \cD_{b_i}) \cong (Y_i, D_i)$ for $i=1,2$, we say that $(Y_1, D_1)$ and $(Y_2, D_2)$ are \emph{deformation-equivalent}.
		\item Let $\mathcal{U}=\cY \setminus \cD$, and let $\pi|_\mathcal{U}\colon \mathcal{U} \rightarrow B$ be the induced deformation of log Calabi--Yau varieties. Given a path $\gamma: [0,1] \rightarrow B$ in the analytic topology with $\gamma(0)=b_1$ and $\gamma(1)=b_2$, parallel transport in the local systems $(R^k\pi|_{\mathcal{U}})_*\ZZ$ induces isomorphisms
		\[
		\gamma_* \colon H^k(U_{b_1}, \ZZ) \cong H^k(U_{b_2}, \ZZ).
		\]
		compatible with cup product, 
		which we will call \emph{parallel transport}.
	\end{itemize}
\end{definition}
\subsection{Looijenga pairs}\label{sec:LP}
In this section, we focus on the two-dimensional case.\\
Following~\cite{GrossHackingKeel} or \cite{Friedman}, a \emph{Looijenga pair} is a smooth projective surface $Y$ together with a connected, singular, nodal anticanonical divisor $D$. So $D$ is either an irreducible nodal curve, or a cycle $D=D_1+\dots+D_n$ of $n$ rational curves. $Y$ is necessarily rational, so we have an isomorphism $\Pic(Y) \cong H^2(Y, \ZZ)$.
It follows from the classification of surfaces that after a sequence of toroidal blowups, any Looijenga pair admits a toric model (\cite[Proposition 1.3]{GrossHackingKeelMS}).
Define the lattice 
\[
\langle D_1, \dots D_n\rangle ^\perp=\{L \in \Pic(Y) \mid L \cdot D_i=0 \; \text{for all}\; i \}
\]
A cyclic ordering of the components of $D$ induces a canonical identification $\Pic^0(D) \cong \CC^\times$ by \cite[Lemma 2.1]{GrossHackingKeel}. We introduce the period point and the marked period point:
\begin{definition} Let $(Y, D)$ be a Looijenga pair.
	\begin{itemize}
		\item The map
		\[
		\phi_Y \colon \langle D_1, \dots D_n\rangle ^\perp \rightarrow \Pic^0(D) \cong \CC^\times, \quad L \mapsto L_{|D}
		\]
		is called the {\it period point} $\phi_Y \in \Hom(\langle D_1, \dots D_n\rangle ^\perp, \CC^\times)$ of $(Y, D)$.
		\item A \emph{marking} of $D$ is a choice of point $p_i \in D_i^\circ$ for every $i$.
		\item	The \emph{marked period point} of $D$ is 
		\[
		\phi_{Y, p_i} \colon \Pic(Y) \rightarrow \Pic^0(D) \cong \CC^\times, \quad L \mapsto L|_{D}^{-1} \otimes \bigotimes_{i=1}^n\cO_{D}((D_i \cdot L)p_i).
		\]
	\end{itemize}
\end{definition}
We define the generic nef cone (c.f \cite[Definition 1.7]{GrossHackingKeel} and \cite[Definition 5.4]{Friedman}).
\begin{definition}
	Let $(Y, D)$ be a Looijenga pair. Let $C^+$ be the positive cone, that is, the connected component of the cone $\{ x \in \Pic(Y )_\RR \mid x^2>0\}$ containing the ample classes.
	For a given ample divisor $H$, we denote $\tilde{\mathcal{M}} \subset \Pic(Y)$ the set of classes $E$ such that 
	$E^2 = K_Y \cdot E = -1$, and $E \cdot H> 0.$\footnote{By \cite[Lemma 2.13]{GrossHackingKeel}, the set $\tilde{\mathcal{M}}$ is independent of $H$.}
	The generic nef cone of $Y$ is 
	\[
	\overline{\mathcal{A}}^{gen}=\left \{ x \in {C}^+ \mid x \cdot D_i \geq 0 \;\text{and} \;x \cdot E \geq 0 \;\text{for all} \;E \in\tilde{\mathcal{M}} \right \}.
	\]
\end{definition}
Note that $\overline{\mathcal{A}}^{gen}$ is denoted $C_D^{++}$ in \cite{GrossHackingKeel}. By \cite[Lemma 2.15]{GrossHackingKeel}, the nef cone $\Nef(Y)$ is the subcone of $\overline{\mathcal{A}}^{gen}$ defined by the inequalities $x \cdot [C] \geq 0$ for all $(-2)$-curves $C$.
In fact, $\overline{\mathcal{A}}^{gen}$ is identified with the nef cone of a very general deformation of the pair $(Y, D)$, see the discussion in \cite{Friedman} above Lemma 5.5.
We recall the global Torelli Theorem of Gross--Hacking--Keel, in the formulation of \cite[Theorem~8.7]{Friedman}.
\begin{theorem}\cite[Theorem 1.8]{GrossHackingKeel} \label{thm:GHKTorelli}
	Let $(Y,D)$ and $(Y',D')$ be Looijenga pairs and let 
	\[
	\mu \colon \Pic(Y) \rightarrow \Pic(Y')
	\] an isomorphism of lattices. 
	Suppose that 
	\begin{enumerate}
		\item	$\mu([D_i])=[D_i']$
		\item $\mu(\overline{\mathcal{A}}^{gen}(Y))=\overline{\mathcal{A}}^{gen}(Y')$
		\item	$\phi_{Y'} \circ \mu = \phi_{Y}$.
	\end{enumerate}
	Then there exists an isomorphism of pairs $f:(Y,D) \to (Y',D')$.
\end{theorem}
To conclude that $f$ induces $\mu$ in Theorem~\ref{thm:GHKTorelli}, one must add the condition that $\mu(\Nef(Y))=\Nef(Y')$.
\subsection{Extensions of mixed Hodge structures}
We write $\ZZ(1)$ for the \emph{Tate Hodge structure}. Up to isomorphism, this is the unique pure Hodge structure on $\ZZ$ of weight $-2$. Given a mixed Hodge structure (MHS) $V$, we write $V(n)$ for $V \otimes \ZZ(1)^{\otimes n}$, this has the effect of raising all weights by $-2n$.
Suppose that 
\begin{equation}\label{eq:ext3}
	0 \rightarrow A \xrightarrow{i} E \xrightarrow{p} B \rightarrow 0
\end{equation}
is an extension in the category of MHS, with $A$ and $B$ pure of weights $a$ and $b$. Suppose that $a <b$. By \cite[Proposition 2]{Carlson}, we have an isomorphism 
\[
\Ext^1_{\MHS}(B, A) \cong J^0(\Hom(B, A))
\] 
where $J^0$ denotes a Jacobian variety. Under our assumptions on the weights we have
\[
J^0(\Hom(B, A)) \cong \frac{\Hom(B, A)_\CC}{\Hom(B, A)_\ZZ},
\]
so any extension of the form \eqref{eq:ext3} has an associated \emph{extension class} $\psi \colon B \rightarrow A \otimes \CC^\times$. This isomorphism can be made explicit: choose a section $s_F \colon B_\CC \rightarrow E_\CC$ of $p$ preserving the Hodge filtration, and a retraction $r_\ZZ \colon E_\ZZ \rightarrow A_\ZZ$ of $i$. 
Then the image of $r_\ZZ \circ s_F$ in 
$J^0(\Hom(B, A))$
is independent of the choices made, and is the corresponding extension class (\cite[Lemma 4]{Carlson}).

If $(Y, D)$ is a Looijenga pair, then the excision long exact sequence gives an exact sequence
\[
0 \rightarrow \ZZ(2) \rightarrow H_2(U, \ZZ) \rightarrow \langle D_1, \dots D_n\rangle ^\perp(2) \rightarrow 0
\]
where $U=Y \setminus D$, so that $H_2(U)$ is an extension of a weight $-2$ Hodge structure by a Hodge structure of weight $-4$. We have the following result.
\begin{proposition}\cite[Proposition 3.12]{Friedman}\label{prop:LooijengaPeriodpoint}
	The extension class associated to $H_2(U)$ is the period point $\phi_Y$.
\end{proposition}
We easily deduce from this that the Torelli theorem as formulated in Conjecture~\ref{conj:torelli1} holds for Looijenga pairs:
\begin{theorem} \label{thm:TorelliLooijengaPairs}
	Conjecture~\ref{conj:torelli1} holds for $(\cY, \cD) \rightarrow B$ a family of Looijenga pairs.
\end{theorem}
\begin{proof}
	For notational simplicity, we write $Y_i=\cY_{b_i}, D_i=\cD_{b_i}$, and $U_i=Y_i \setminus D_i$.
	Let 
	\[
	\gamma_* \colon H_2(U_1) \rightarrow H_2(U_2)
	\]
	be the isomorphism of MHS induced by parallel transport.
	Since the deformation $\cD \rightarrow B$ preserves the combinatorial type of the dual complex, $\gamma_*([D_i])=[D_i']$. The condition that $\gamma_*$ is an isomorphism of MHS implies that 
	\[
	\psi_{Y_2} \circ \gamma_* = \psi_{Y_1}
	\]
	where  $\psi_{Y_i}$ is the extension class of $H_2(U_i)$. By Proposition~\ref{prop:LooijengaPeriodpoint}, the extension classes are identified with the period points $\phi_{Y_i}$, so we obtain $\phi_{Y_2} \circ \gamma_* = \phi_{Y_1}$.\\
	Finally, \cite[Lemma 2.13]{GrossHackingKeel} shows that the cone $\overline{\mathcal{A}}^{gen}(Y)$ (denoted $C^{++}_D$ there) is preserved by parallel transport, so we conclude by Theorem~\ref{thm:GHKTorelli}.
\end{proof}
\subsection{Log Calabi--Yau threefolds}
Suppose now that $(Y, D)$ is a log CY pair with $\dim(Y)=3$. By the adjunction formula, we have that 
\[
\omega_D=(\omega_Y \otimes \cO(D))|_D \cong \cO_D.
\] 
Throughout this paper, we will only consider the case that $D$ has maximal boundary, i.e the dual complex has the maximal possible dimension $2$. This places several restrictions on $D$.  
We have the following result
\begin{proposition}\label{lem:logcy3}
	Let $(Y, D)$ be a three-dimensional log Calabi--Yau pair with maximal boundary.
	The dual complex $\Sigma_D$ is homeomorphic to $S^2$. 
	Moreover, each component of $D$ is a rational surface, and the double curves along each component form a cycle of rational curves on $S$. 
\end{proposition}
\begin{proof}
	The fact that the dual complex is homeomorphic to $S^2$ is \cite[33]{KollarXu}. For the remaining claims, consider a component $S$ of $D$. The divisor $B=(D-S)|_S$ is an anticanonical divisor on $S$, and by \cite[Section 2.2]{Kaloghiros}, the pair $(S, B)$ is a log Calabi--Yau surface pair with maximal boundary. By the classification of surfaces, $S$ is a rational surface, and $B$ is a cycle of smooth rational curves (note that $B$ irreducible nodal is ruled out by the assumption that $D$ is snc). 
\end{proof}
In analogy with the Kulikov classification of degenerations of $K3$ surfaces~\cite{Kulikov}, we define 
\begin{definition}
	A \emph{type III surface} is a reduced simple normal crossing surface $D$ with $\omega_D \cong \cO_D$ such that the dual complex is homeomorphic to $S^2$. 
\end{definition}
We will study type III surfaces in detail and will use the following notation:
\begin{definition}
	Let $D$ be a type III surface with dual complex $\Sigma$.
	\begin{itemize}
		\item We denote $\Sigma^{[i]}$ the set of strata of $\Sigma$ of (real) dimension $i$, and write $\Sigma(i)$ for the cardinality of the set $|\Sigma^{[i]}|$. 
		\item We write $D_v$ for the irreducible component of $D$ corresponding to $v \in \Sigma^{[0]}$, and $D_e$ for the one-dimensional stratum of $D$ corresponding to $e \in \Sigma^{[1]}$.
		\item We recall that $D_e^{\circ}$ denotes $D_e$ with the two zero-dimensional strata removed, so that $D_e^{\circ} \cong \CC^\times$. 
		\item 	Given a component $D_v$, write $\partial D_v$ for the union of one-dimensional strata of $D$ contained in $D_v$. $\partial D_v$ is an anticanonical cycle of rational curves, so that $(D_v, \partial D_v)$ is a Looijenga pair. 
		\item An \emph{orientation} of $D$ is an orientation of $\Sigma$, i.e a choice of generator for $H_2(S^2, \ZZ)$. An orientation of $D$ induces an orientation for each Looijenga pair $(D_v, \partial D_v)$, and therefore canonical isomorphisms 
		$\Pic^0(\partial D_v) \cong \CC^\times$ by \cite[Lemma 2.1]{GrossHackingKeel}. 
		\item An \emph{orientation of edges} is a choice of direction for every one-dimensional stratum of $\Sigma$.
		\item A \emph{marking} of $D$ is a choice of point $p_e \in D_e^{\circ}$ for each $e \in \Sigma^{[1]}$.
		\item Given $e \in \Sigma^{[1]}$, we write $e=e(v,w)$ to mean that $v$ and $w$ are the vertices on $e$. Viewing $D_e$ as a boundary component of either $D_v$ or $D_w$ gives two identifications $D_e^{\circ} \cong \CC^\times$. If $p \in D_e^{\circ}$ corresponds to $\lambda \in \CC^\times$ under one identification, $p$ corresponds to $\lambda^{-1}$ under the other identification, so we unambiguously define $m_e$ to be the point on $D_e$ that corresponds to $-1$ under either identification.
	\end{itemize}
\end{definition}
We recall the two different types of interior blowups from Lemma~\ref{lem:centersoftoricmodel}.
\begin{itemize}\label{interiorblowups}
	\item Let $p \in D_e$ be a point in the interior of a one-dimensional stratum, and let $Y'$ be the blowup of $Y$ in $p$. Denoting by $D'$ the strict transform of $D$, the pair $(Y', D')$ is again log Calabi--Yau. The exceptional divisor $E$ meets the two components of $D'$ containing $D_e'$ in a smooth rational curve of self-intersection $(-1)$, and intersects all other components trivially (see Figure~\ref{fig1}).
	\item Let $C \subset D_{v_0}$ be a curve meeting one-dimensional strata transversely and avoiding zero-dimensional strata, and let $Y'$ be the blowup of $Y$ in $C$. Denoting by $D'$ the strict transform of $D$, the pair $(Y', D')$ is again log Calabi--Yau. The exceptional divisor $E$ meets $D_{v_0} \cong D'_{v_0}$ in a section of $E \rightarrow C$, and adjacent components $D'_v$ in a (possibly empty) set of disjoint smooth rational curves $E^i_v$, which are $(-1)$-curves on $D'_v$ (see Figure~\ref{fig2}).
\end{itemize}       
\begin{figure}
	\tikzset{every picture/.style={line width=0.75pt}} 
	\tikzset{every picture/.style={line width=0.75pt}} 
	
	\begin{tikzpicture}[x=0.75pt,y=0.75pt,yscale=-1,xscale=1]
		
		\draw [line width=1.5]    (589.5,66) -- (547,110) ;
		\draw [line width=1.5]    (547,110) -- (548.5,208) ;
		\draw [line width=1.5]    (507.5,68) -- (547,110) ;
		\draw    (338.5,110) -- (428,109.02) ;
		\draw [shift={(430,109)}, rotate = 179.37] [color={rgb, 255:red, 0; green, 0; blue, 0 }  ][line width=0.75]    (10.93,-3.29) .. controls (6.95,-1.4) and (3.31,-0.3) .. (0,0) .. controls (3.31,0.3) and (6.95,1.4) .. (10.93,3.29)   ;
		\draw [line width=1.5]    (548.5,208) -- (588,250) ;
		\draw [line width=1.5]    (548.5,208) -- (506,252) ;
		\draw  [fill={rgb, 255:red, 144; green, 19; blue, 254 }  ,fill opacity=1 ] (544,172.5) .. controls (544,170.57) and (545.68,169) .. (547.75,169) .. controls (549.82,169) and (551.5,170.57) .. (551.5,172.5) .. controls (551.5,174.43) and (549.82,176) .. (547.75,176) .. controls (545.68,176) and (544,174.43) .. (544,172.5) -- cycle ;
		\draw [line width=1.5]    (257.5,68) -- (215,112) ;
		\draw [line width=1.5]    (215,112) -- (216.5,210) ;
		\draw [line width=1.5]    (175.5,70) -- (215,112) ;
		\draw [line width=1.5]    (216.5,210) -- (256,252) ;
		\draw [line width=1.5]    (216.5,210) -- (174,254) ;
		\draw  [fill={rgb, 255:red, 144; green, 19; blue, 254 }  ,fill opacity=1 ] (212,173.5) .. controls (212,171.57) and (213.68,170) .. (215.75,170) .. controls (217.82,170) and (219.5,171.57) .. (219.5,173.5) .. controls (219.5,175.43) and (217.82,177) .. (215.75,177) .. controls (213.68,177) and (212,175.43) .. (212,173.5) -- cycle ;
		\draw  (177.34,152.23) .. controls (184.6,137.67) and (199.08,127.75) .. (215.75,127.75) .. controls (231.77,127.75) and (245.76,136.91) .. (253.27,150.53) -- (215.75,173.5) -- cycle ;
		\draw [color={rgb, 255:red, 144; green, 19; blue, 254 }  ,draw opacity=1 ][line width=2.25]    (177.34,152.23) -- (215.75,173.5) ;
		\draw [color={rgb, 255:red, 144; green, 19; blue, 254 }  ,draw opacity=1 ][line width=2.25]    (215.75,173.5) -- (253.27,150.53) ;
		
		\draw (238.51,159.01) node [anchor=north west][inner sep=0.75pt]  [font=\footnotesize] [align=left] {$\displaystyle -1$};
		\draw (172,162) node [anchor=north west][inner sep=0.75pt]  [font=\footnotesize] [align=left] {$\displaystyle -1$};

	\end{tikzpicture}
	
	\caption{Two components of $D$ meeting along a one-dimensional stratum $D_e$ (the vertical line). The picture shows an interior blowup of a point on $D_e$.}
	\label{fig1}
\end{figure}
\begin{figure}

	\tikzset{every picture/.style={line width=0.75pt}} 
	
	\begin{tikzpicture}[x=0.75pt,y=0.75pt,yscale=-1,xscale=1]
		
		\draw [line width=1.5]    (178,28) -- (179,81) ;
		\draw [color={rgb, 255:red, 144; green, 19; blue, 254 }  ,draw opacity=1 ][line width=2.25]    (126,157) -- (155.25,174.76) -- (182,191) ;
		\draw [color={rgb, 255:red, 144; green, 19; blue, 254 }  ,draw opacity=1 ][line width=2.25]    (182,191) -- (201,230) ;
		\draw [color={rgb, 255:red, 144; green, 19; blue, 254 }  ,draw opacity=1 ][line width=2.25]    (225.5,158) -- (193,182) ;
		\draw [color={rgb, 255:red, 144; green, 19; blue, 254 }  ,draw opacity=1 ][line width=2.25]    (126,124) -- (126,142) -- (126,157) ;
		\draw    (182,158) -- (182,191) ;
		\draw [color={rgb, 255:red, 144; green, 19; blue, 254 }  ,draw opacity=1 ][line width=2.25]    (225.5,125) -- (225.5,158) ;
		\draw [color={rgb, 255:red, 144; green, 19; blue, 254 }  ,draw opacity=1 ][line width=2.25]    (229.5,257) -- (201,230) ;
		\draw    (182,158) -- (229.5,257) ;
		\draw    (225.5,125) -- (182,158) ;
		\draw    (126,124) -- (136.84,130.61) -- (182,158) ;
		\draw [line width=1.5]    (179,81) -- (142,132) ;
		\draw [line width=1.5]    (77,225) -- (272,235) ;
		\draw [line width=1.5]    (179,81) -- (213,136) ;
		\draw [line width=1.5]    (225.5,158) -- (272,235) ;
		\draw [line width=1.5]    (126,157) -- (77,225) ;
		\draw [line width=1.5]    (272,235) -- (321,261) ;
		\draw [line width=1.5]    (30,253) -- (77,225) ;
		\draw [line width=1.5]    (503,27) -- (504,80) ;
		\draw [color={rgb, 255:red, 144; green, 19; blue, 254 }  ,draw opacity=1 ][line width=1.5]    (451,156) -- (507,190) ;
		\draw [color={rgb, 255:red, 144; green, 19; blue, 254 }  ,draw opacity=1 ][line width=1.5]    (507,190) -- (526,229) ;
		\draw [color={rgb, 255:red, 144; green, 19; blue, 254 }  ,draw opacity=1 ][line width=1.5]    (550.5,157) -- (507,190) ;
		\draw [line width=1.5]    (504,80) -- (451,156) ;
		\draw [line width=1.5]    (402,224) -- (597,234) ;
		\draw [line width=1.5]    (504,80) -- (550.5,157) ;
		\draw [line width=1.5]    (550.5,157) -- (597,234) ;
		\draw [line width=1.5]    (451,156) -- (402,224) ;
		\draw [line width=1.5]    (597,234) -- (646,260) ;
		\draw [line width=1.5]    (355,252) -- (402,224) ;
		\draw (519.5,185) node [anchor=north west][inner sep=0.75pt]  [font=\footnotesize] [align=left] {$\displaystyle \textcolor[rgb]{0.56,0.07,1}{C \cong \PP^1 }$};
		\draw (103,132) node [anchor=north west][inner sep=0.75pt]  [font=\footnotesize] [align=left] {$\displaystyle -1$};
		\draw (193,242) node [anchor=north west][inner sep=0.75pt]  [font=\footnotesize] [align=left] {$\displaystyle -1$};
		\draw (229,135) node [anchor=north west][inner sep=0.75pt]  [font=\footnotesize] [align=left] {$\displaystyle -1$};
		\draw    (288.5,110) -- (378,109.02) ;
		\draw [shift={(380,109)}, rotate = 179.37] [color={rgb, 255:red, 0; green, 0; blue, 0 }  ][line width=0.75]    (10.93,-3.29) .. controls (6.95,-1.4) and (3.31,-0.3) .. (0,0) .. controls (3.31,0.3) and (6.95,1.4) .. (10.93,3.29)   ;
	\end{tikzpicture}
	
	\caption{One component $D_{v_0}$ of $D$ is the triangle, meeting adjacent components along its edges. The picture shows an interior blowup of a smooth curve $C$ contained in $D_{v_0}$. Note that $C$ is drawn tropically, and represents a smooth $\PP^1$. Similar for the smooth exceptional divisor, which is a $\PP^1$-bundle over $C$.}
	\label{fig2}
\end{figure}
\subsection{Mori theory}
The Mori cone $\overline{\NE}(Y) \subset N_1(Y, \RR)$ is the closed convex cone defined as the closure of the cone spanned by the classes of curves on $Y$, and the nef cone $\Nef(Y) \subset N^1(Y, \RR)$ is its dual cone.
The effective cone $\Eff(Y) \subset N^1(Y, \RR)$ is the convex cone spanned by the classes of effective divisors. Note that $\Eff(Y)$ is usually not closed (c.f \cite[Example 1.5.1]{Lazarsfeld})\\
The varieties considered in this paper will usually be rational and therefore satisfy $h^1(Y, \cO_Y)=h^2(Y, \cO_Y)=0$. Hence, we may equivalently view these cones inside $H_2(Y, \RR)$ or $H^2(Y, \RR)$.
Given an extremal ray $R$ of $\overline{\NE}(Y)$ with $K_Y \cdot R< 0$, the cone theorem shows that there exists a contraction morphism $\pi \colon Y \rightarrow {Y}'$ of relative Picard rank $1$ with connected fibers and satisfying 
$\pi(C)=pt$ for a curve $C$ if and only if $[C] \in R$. For any curve $C$ with $[C]\in R$, \cite[Theorem 3.2]{KollarMori} gives an exact sequence
\begin{equation}\label{eq24}
	0 \rightarrow \Pic({Y}') \xrightarrow{\pi^*} \Pic(Y) \xrightarrow{ \cdot C} \ZZ
\end{equation}
If $\pi \colon (Y, D) \rightarrow (Y', D')$ is an interior blowup of log CY3 pairs, then $\pi$ corresponds to the contraction of an $K_{Y}$-negative extremal ray of $\overline{\NE}(Y)$, since $D \in |-K_{Y}|$.
We have the following well-known result:
\begin{lemma}\label{lem:conespullback}
	Let $\pi \colon Y \rightarrow {Y'}$ be the contraction of a $K_Y$-negative extremal ray $\ell$. 
	\begin{enumerate}[(a)]
		\item $\pi^*$ identifies $\Nef({Y}')$ with the face $\ell^\perp \cap \Nef(Y)$ of $\Nef(Y)$.
		\item $\pi^*$ identifies $\Eff({Y}')$ with $\ell^\perp \cap \Eff(Y) \subset \Eff(Y)$.
	\end{enumerate}
\end{lemma}
\begin{proof}
	The morphism $\pi$ is surjective, so by \cite[Example~1.4.4]{Lazarsfeld}, a class $x \in \Nef({Y}')$ if and only if $\pi^*x \in \Nef(Y)$. Using \eqref{eq24}, this proves a).\\
	For b), it is clear that if $x \in \Eff({Y}')$ then $\pi^*x \in \Eff(Y)$. Conversely, if $\pi^*x \in \Eff(Y)$, then $\pi_*\pi^*x \in \Eff({Y}')$, and $\pi_*\pi^*x=x$ since $\pi$ is surjective. We conclude again using \eqref{eq24}.
\end{proof}
Suppose now that $\pi \colon Y \rightarrow {Y}'$ is the contraction of a $K_Y$-negative extremal ray on a smooth threefold, and let $\ell$ be the first lattice point on $R$ represented by an effective curve. 
Suppose further that $\pi$ is a divisorial contraction, and denote $E$ the exceptional divisor. By Mori's classification of contractions of $K_Y$-negative extremal rays~\cite[Theorem 3.3]{Mori}, the contraction $\pi$ is one of the following 
\begin{enumerate}[(1)]\label{Moricontractions}
	\item $\pi$ is the blowup of a smooth curve $C$, $E$ is the exceptional divisor (a $\PP^1$-bundle over $C$), and $\ell$ is the class of the ruling.
	\item $\pi$ is the blowup of a smooth point, $E \cong \PP^2$ is the exceptional divisor, $\cO_E(E) \cong \cO_{\PP^2}(-1)$, and $\ell$ is the class of a line in $\PP^2$.
	\item $\pi$ is the blowup of a singularity of type $A_1$ (an ordinary double point), $E \cong \PP^1 \times \PP^1$ is the exceptional divisor, $\cO_E(E)$ has bidegree $(-1,-1)$, and $\ell$ is the class of either ruling.
	\item $\pi$ is the blowup of a singularity of type $A_2$, $E$ is the exceptional divisor, $E$ is isomorphic to a singular quadric in $\PP^3$, $\cO_E(E)=\cO_{\PP^3}(-1)|_E$, and $\ell$ is the class of a ruling of the quadric.
	\item $\pi$ contracts $E=\PP^2$ to a singular point, $\cO_E(E)=\cO_{\PP^2}(-2)$, and $\ell$ is the class of a line in $\PP^2$.
\end{enumerate}
We easily compute the following intersection numbers
\begin{figure}
	\begin{equation*}
		\begin{tabular}{ c|c|c|c}
			Type & $E \cdot \ell$ & $K_Y \cdot \ell$&$E \cdot K_Y^2$ \\ 	
			\hline
			$1$ & $-1$ &$-1$ &$a$\\ 
			$2$ & $-1$  & $-2$ & 4\\ 
			$3$ & $-1$  & $-1$ & 2\\ 
			$4$ & $-1$  & $-1$ &2\\ 
			$5$ & $-2$ & $-1$&1 \\ 
		\end{tabular}
	\end{equation*}
	\caption{Numerical properties of the birational Mori contractions.}
	\label{table}
\end{figure}
where for type $1$, the integer $a$ depends on the genus of $C$ and the degree of the normal bundle $N_{C/{Y'}}$.
It follows from the negativity lemma~\cite[Lemma 3.39]{KollarMori} that the exceptional divisor $E$ is irreducible. The class $[E]$ therefore spans an extremal ray of $\Eff(Y)$. We have the following result:
\begin{lemma}\label{lem:primitive}
	Let $Y$ be a smooth projective threefold, and let $\pi \colon Y \rightarrow {Y}'$ be the divisorial contraction of a $K_Y$-negative extremal ray $R$ of $\overline{\NE}(Y)$. Let $E$ be the exceptional divisor, and let $\ell$ be the first lattice point on $R$ represented by an effective curve. Then 
	\begin{enumerate}[(a)]
		\item $\ell$ is the primitive generator of the extremal ray $R$ of $\overline{\NE}(Y)$. 
		\item $[E]$ is the primitive generator of the extremal ray $\RR_{\geq 0}[E]$ of $\Eff(Y)$. 
	\end{enumerate}
\end{lemma}
\begin{proof}
	For a), it suffices to find a divisor $D$ with $D \cdot \ell=\pm 1$. By the table in Figure~\ref{table}, 
	we may take $D=E$ in cases $1-4$, and $D=K_Y$ in case $5$.
	For part b), we need to find a curve $C$ with $E \cdot C=\pm 1$. In cases $1-4$, we may take $C=\ell$. 
	In case $5$, we see from the table in Figure~\ref{table} that $E \cdot K_Y^2=1$, so $E$ is not divisible, as required. 
\end{proof}
\section{Torelli for type III  surfaces}
Let $D$ be a type III surface. Throughout we fix an orientation of $D$, as well as an orientation of edges. Consider the homomorphism
\[
\ell \colon 
\bigoplus_{v \in \Sigma^{[0]}} \Pic(D_v) \rightarrow \bigoplus_{e \in \Sigma^{[1]}} \Pic(D_e), \qquad (L_v)_{v \in \Sigma^{[0]}} \mapsto (L_v|_{D_{e}}-L_w|_{D_{e}})_{e=e(v,w) \in \Sigma^{[1]}}
\]
where the signs are determined by the orientation of edges.
The MHS on $H^2(D)$ was studied in detail by Friedman--Scattone~\cite[p.16]{FriedmanScattone}. It has weights $0$ and $2$ and fits in an extension 
\begin{equation}\label{eq:ext4}
	0 \rightarrow \ZZ(0) \xrightarrow{i} H^2(D) \xrightarrow{p} \Lambda \rightarrow 0.
\end{equation}
where $\Lambda=\ker(\ell)$.
The corresponding extension class is a homomorphism $\psi \colon \Lambda \rightarrow \CC^\times$. We will show that just as in the two-dimensional case, $\psi$ has a geometric description.
\subsection{The period point of $D$}
\begin{definition}
	Let $D$ be a type III surface with marking $p_e$. The marked period point of $D$ is
	\[
	\phi_{D, p_e} \colon \bigoplus_{v \in \Sigma^{[0]}} \Pic(D_v) \rightarrow \CC^\times, \quad (L_v)_{v \in \Sigma^{[0]}} \mapsto \prod_{v \in \Sigma^{[0]}}\phi_{D_v, p_e}(L_v)
	\]
	where $\phi_{D_v, p_e}$ is the marked period point of the Looijenga pair $(D_v, \partial D_v)$.
\end{definition}
Upon restriction to $\Lambda$, the period point is independent of the marking:
\begin{lemma}\label{lem:periodpointindependent}
	The restriction of $\phi_{D, p_e}$ to $\Lambda$ is independent of the choice of marking $p_e$.
\end{lemma}
\begin{proof}
	Consider an edge $e=e(v,w)$. If $p_e$ corresponds to $\lambda \in \CC^\times$ under the identification $\CC^\times \cong \Pic^0(\partial D_v)$, then $p_e$ corresponds to $\lambda^{-1}$ under the identification $\CC^\times \cong \Pic^0(\partial D_w)$. If $l:=(L_u)_{u \in \Sigma^{[0]}} \in \Lambda$, then $L_u \cdot D_e=0$ for $u \neq v$ or $w$, and $L_v \cdot D_e=L_w \cdot D_e$. Therefore, the contribution of $p_e$ in $\phi_{D, p_e}(l)$ cancels.
\end{proof}
Given Lemma~\ref{lem:periodpointindependent}, we define 
\begin{definition}
	Let $D$ be a type III surface. The (unmarked) period point of $D$ is the homomorphism
	\[
	\phi_D \colon \Lambda \rightarrow \CC^\times
	\]
	obtained by restricting $\phi_{D, p_e}$ to $\Lambda$, for an arbitrary choice of marking $p_e$.
\end{definition}
Just as for Looijenga pairs, the period point agrees with the extension class of the corresponding mixed Hodge structure. The following argument is taken from (\cite{Siebert}, private communication)
\begin{proposition}\label{prop:periodpointsagree}
	The extension class $\psi \colon \Lambda \rightarrow \CC^\times$ of the MHS $H^2(D)$ agrees with the (unmarked) period point $\phi_D$.
\end{proposition}
\begin{proof}
	Let $D^{[i]}$ denote the disjoint union of strata of $D$ of codimension $i$. We have a resolution of $\cO_{D}^*$
	\[
	1 \rightarrow \cO_D^* \rightarrow \cO_{D^{[0]}}^* \rightarrow \cO_{D^{[1]}}^* \rightarrow \cO_{D^{[2]}}^*\rightarrow 1
	\]
	(where the signs in the differentials are determined by fixing an ordering of the vertices of $\Sigma_D$).
	The hypercohomology spectral sequence associated to this resolution has $E_1$-page
	\[
	\begin{tikzpicture}
		\matrix (m) [matrix of math nodes,
		nodes in empty cells,nodes={minimum width=5ex,
			minimum height=5ex,outer sep=-5pt},
		column sep=1ex,row sep=1ex]{
			&      &     &     &\\
			1     &  H^1(\cO^*_{D^{[0]}}) &  H^1(\cO^*_{D^{[1]}}) & 0 & \\
			0     &  H^0(\cO^*_{D^{[0]}})  & H^0(\cO^*_{D^{[1]}})&  H^0(\cO^*_{D^{[2]}})\\
			\quad\strut &   0  &  1  &  2  &  \strut \\};
		\draw[thick] (-4,-1.25) -- (4,-1.25);
		\draw[thick] (-3.25,-2) -- (-3.25,2);
	\end{tikzpicture}
	\]
	and we have that 
	\[
	\ker\left(H^1(\cO^*_{D^{[0]}}) \rightarrow H^1(\cO^*_{D^{[1]}}) \right) \cong \ker\left(\Pic(D^{[0]}, \ZZ) \rightarrow \Pic(D^{[1]}, \ZZ)\right)=\Lambda
	\]
	so that the differential $d_2$ induces a map 
	\[
	d_2 \colon \Lambda \rightarrow H^2(|\Sigma|, \CC^\times) \cong \CC^\times.
	\]
	By \cite[Proposition 3.4]{FriedmanScattone}, $d_2$ is identified with the extension class $\psi$. 
	We can describe the differential $d_2$ explicitly. Let $l = (L_v)_{v \in \Sigma^{[0]}} \in \Lambda$. For each edge $e=e(v,w)$, we have $L_v \cdot D_e=L_w \cdot D_e$ by definition of $\Lambda$, and since $D_e \cong \PP^1$, we choose for each edge an isomorphism
	\begin{equation}\label{eq:iso2}
		\psi_{e} \colon L_v|_{D_e} \simeq L_w|_{D_e}.
	\end{equation}
	For each $\sigma \in \Sigma^{[2]}$, let $p_\sigma \in D$ be the corresponding zero-dimensional stratum of $D$.  
	Let $e_1, e_2, e_3$ be the edges of $\sigma$ in the cyclic order given by the orientation of $\sigma$ induced by the orientation of $\Sigma$, and write $e_i=e_i(v_i, v_{i+1})$ (with indices modulo $3$). Then $\psi_{e_3}|_{p_\sigma}\circ  \psi_{e_2}|_{p_\sigma} \circ \psi_{e_1}|_{p_\sigma}$ is an automorphism of $L_{v_1}|_{p_\sigma}$ given by an nonzero scalar $\alpha_\sigma \in \CC^\times$ (independent of the choice of $v_1$). Then 
	\begin{equation}\label{eq1}
		d_2(l)=\prod_{\sigma \in \Sigma^{[2]}}\alpha_\sigma  \in \CC^\times
	\end{equation}
	Fix a marking $p_e$ of $D$. For $e=e(v,w)$ write
	\[
	d_e =L_v \cdot D_e =L_w \cdot D_e \in \ZZ.
	\]
	We recall that 
	$\phi_{D_v, p_e} (L_v)$ can be computed explicitly as follows: choose for each component $D_e$ of $\partial D_v$ a nowhere vanishing section $s_{v,e}$ of $L_v^{-1}|_{D_e} \otimes \cO_{D_e}(d_e p_e)$, and set 
	\[
	\lambda_e=\frac{s_{v,e}(0)}{s_{v,e}(\infty)}
	\]
	where $0$ and $\infty$ are the points on $D_e$ corresponding to the zero-dimensional strata, with the order determined by the orientation. 
	Then by \cite[Lemma 2.2 (1)]{GrossHackingKeel}, we have $\phi_{D_v}(L_v)=\prod_{e \colon v \in e} \lambda_e$. 
	We may rewrite $\phi_{D}(l)=\prod_{v \in \Sigma^{[0]}}\phi_{D_v}(L_v)$ as a product over $\sigma \in \Sigma^{[2]}$ of terms of the form 
	\begin{equation}\label{eq:ratio}
		\frac{s_{v,e}(p_\sigma)}{s_{w,e}(p_{\sigma})}\frac{s_{w,e}(p_{\sigma})}{s_{u,e}(p_{\sigma})}\frac{s_{u,e}(p_{\sigma})}{s_{v,e}(p_{\sigma})}
	\end{equation}
	The ratio $\frac{s_{v,e}}{s_{w,e}}$ defines an isomorphism $L_v|_{D_e} \cong L_w|_{D_e}$, so that the scalar \eqref{eq:ratio} is equal to $\alpha_\sigma$. 
	We therefore obtain $\psi(l)=\phi_{D}(l)$, as required.
\end{proof}
We have the following Corollary:
\begin{corollary}\label{cor:periodpoint}
	Let $\phi_D \colon \Lambda \rightarrow \CC^\times$ be the (unmarked) period point of a type III surface $D$. 
	\begin{enumerate}
		\item $\ker(\phi_D)=\Pic(D)$.
		\item If for each $v \in \Sigma^{[0]}$, the Looijenga pair $(D_v, \partial D_v)$ is toric, then $\phi_D$ is the trivial homomorphism, in particular $\Pic(D)=\Lambda$.
	\end{enumerate}
\end{corollary}
\begin{proof}
	By Theorem~\ref{prop:periodpointsagree}, $\phi_D$ is identified with the extension class $\psi$ of the mixed Hodge structure $H^2(D)$, so (1) is \cite[Proposition 3.1]{FriedmanScattone}. For $(2)$,
	choose as marking the points $m_e$ that correspond to $-1$ on each one-dimensional stratum $D_e$.
	Since each Looijenga pair $(D_v, \partial D_v)$ is toric, \cite[Lemma 2.8]{GrossHackingKeel} shows 
	\[
	L_v|_{{D}_e} \cong \bigotimes_{e \colon v \in e} \cO_{\partial {D}_v}((L_v \cdot D_e)m_e)
	\]
	and it follows that the marked period point $\phi_{{D}_v, m_e} \colon \Pic({D}_v) \rightarrow \CC^\times$ is trivial. 
	It follows that $\phi_D \colon \Lambda \rightarrow \CC^\times$ is the trivial homomorphism.
\end{proof}
\subsection{Torelli for type III surfaces}
Given our identification of the period point with the extension class of $H^2(D)$, we can prove a Torelli theorem for type III surfaces. We start with a preparatory lemma. Recall that $\Aut^0$ denotes the identity component of the automorphism group.\\
Fix a component $D_v$, an element $\alpha \in \Aut^0(\partial D_v)$, and $L \in \Pic(\partial D_v)$.  Following \cite[Lemma 2.5]{GrossHackingKeel} we define 
a homomorphism $\theta_v \colon \Aut^0(\partial D_v) \rightarrow \Hom(\Pic(D_v), \CC^\times)$
via
\[
\theta_v(\alpha)(L) = L|_{\partial D_v}^{-1} \otimes \alpha^*L|_{\partial D_v} \in \Pic^0(\partial D_v) \cong \CC^\times.
\]
Denote $D^1$ the $1$-skeleton of $D$ (note that $D^1$ is connected, as opposed to $D^{[1]}$). An element $\alpha \in \Aut^0(D^1)$ is an automorphism of $D^1$ fixing zero-dimensional strata pointwise, so restricts to an element $\alpha|_{\partial D_v} \in \Aut^0(\partial D_v)$ for each $v \in \Sigma^{[0]}$. The homomorphisms $\theta_v$ assemble to a homomorphism
\[
\theta \colon \Aut^0(D^1) \rightarrow \Hom(\Pic(D^{[0]}), \CC^\times)
\]
defined as 
\[
\theta(\alpha)(l)=\prod_{v \in \Sigma^{[0]}}\theta_v(\alpha|_{\partial D_v})(L_v)
\]
for $l=(L_v)_{v \in \Sigma^{[0]}}$.
We have the following result:
\begin{lemma}
	Let $D$ be a type III surface, and denote $D^1$ the $1$-skeleton of $D$.
	The orientation of edges of $D$ induces a canonical isomorphism $\Aut^0(D^1) \cong (\CC^\times)^{\Sigma(1)}$. Under this identification, $\theta$ becomes:
	\begin{equation}\label{eq3}
		\theta((\lambda_e)_{e \in \Sigma^{[1]}})(l)=\prod_{e=e(v,w) \in \Sigma^{[1]}} \lambda_e^{(L_v-L_w) \cdot D_e}
	\end{equation}
	and we have an exact sequence
	\begin{equation}\label{eq:seq2}
		\Aut^0(D^1) \xrightarrow{\theta} \Hom(\Pic(D^{[0]}), \CC^\times) \rightarrow \Hom(\Lambda, \CC^\times) \rightarrow 0
	\end{equation}
\end{lemma}
\begin{proof}
	An element $\Aut^0(D^1)$ is an automorphism of $D^1$ fixing zero-dimensional strata pointwise. Let $D_e \cong \PP^1$ be a component of $D^1$. The orientation of $e$ gives a canonical identification of $D_e^\circ$ with $\CC^\times$, so that the automorphism group of $D_e^\circ$ is canonically identified with $\CC^\times$. \\
	The second part is almost identical to the proof of \cite[Lemma 2.5]{GrossHackingKeel}, with the signs determined by the orientation of the edges: to obtain \eqref{eq:seq2}, we apply the functor $\Hom(-, \CC^\times)$ to the exact sequence 
	\begin{equation}\label{eq:seq1}
		0 \rightarrow \Lambda \rightarrow \Pic(D^{[0]}) \rightarrow \Pic(D^{[1]})
	\end{equation}
	defining $\Lambda$, and identify $\Hom(\Pic(D^{[1]}), \CC^\times)$ with $\Aut^0(D^1)$. Note that exactness is preserved since $\CC^\times$ is divisible.
	The fact that the first map coincides with $\theta$ then follows from \eqref{eq3}.
\end{proof}
We note that 
\[
\theta(\alpha) \phi_{D, \{p_e\}}=\phi_{D, \{\alpha(p_e)\}}
\]
as elements of  $\Hom(\Pic(D^{[0]}), \CC^\times)$,
so that the exact sequence $\eqref{eq:seq2}$ in particular expresses the fact that the set of markings for the period point $\phi_D$ is a torsor for $\Aut^0(D^1)$.
\begin{theorem}\label{thm:typeIIITorelli}
	Let $D$ and $D'$ be type III surfaces, and suppose that parallel transport an isomorphism $\mu \colon H^2(D) \rightarrow H^2(D')$ of mixed Hodge structures. Then $D$ and $D'$ are isomorphic.\\ 
	In particular, every deformation type of type III surfaces contains a unique surface with split mixed Hodge structure.
\end{theorem}
\begin{proof}
	Recall that $H^2(D)$ is an extension of the weight $2$ Hodge structure $\Lambda$ by $\ZZ(0)$. 
	Let $\phi_D$ and $\phi_{D'}$ be the unmarked period points. Parallel transport induces isometries $\mu_v \colon \Pic(D_v) \rightarrow \Pic(D_v')$, so that
	\[
	\oplus_v \mu_v \colon \bigoplus_v \Pic(D_v) \rightarrow \bigoplus_v \Pic(D_v')
	\]
	induces an isomorphism $\oplus_v\mu_v|_{\Lambda} \colon \Lambda \rightarrow \Lambda'$. This isomorphism agrees with the isomorphism on weight $2$ pieces induced by $\mu$. Let ${\psi}$ and ${\psi'}$ be the extension classes of $H^2(D)$ and $H^2(D')$. By assumption, $\mu$ induces an isomorphism of mixed Hodge structures, and therefore preserves the period points:
	\[
	{\psi}={\psi}' \circ \oplus_v\mu_v|_{\Lambda}
	\]
	and by Proposition~\ref{prop:periodpointsagree}, $\psi$ and $\psi'$ are identified with the period points $\phi_D$ and $\phi_{D'}$, so we see that $\oplus_v \mu_v$ preserves the period point.
	Choose arbitrary markings $p_e$ and $p_e'$ for $D$ and $D'$. By the exact sequence \eqref{eq:seq2}, we may adjust the marking on $D$ such that we have an equality
	\[
	\phi_{D, p_e}=\phi_{D', p_e'} \circ \oplus_v\mu_v
	\]
	of marked period points.
	Restricting to $\Pic(D_v)$, we obtain 
	\[
	\phi_{D_v, p_e}=\phi_{D_v', p_e'} \circ \mu_v.
	\]
	Since $\mu_v$ is induced by parallel transport, it preserves the generic nef cone (see \cite[Lemma 2.13]{GrossHackingKeel}, denoted $C_{D}^{++}$ there), so the Torelli theorem for marked Looijenga pairs~\cite[Theorem 8.8]{Friedman}, shows that there exists an isomorphism 
	\[
	f_v \colon (D_v, \partial D_v, p_e) \rightarrow (D_v', \partial D_v', p_e')
	\]
	of marked pairs. 
	Given two components $D_v$ and $D_w$ meeting along an edge $D_e$, the isomorphisms $f_v$ and $f_w$ agree on three points on $D_e$ (the two zero-dimensional strata on $D_e$ and the marked point $p_e$). We therefore have $f_v|_{D_e}=f_w|_{D_e}$, so we may glue the isomorphisms $f_v$ to an isomorphism 
	\[
	f \colon D \rightarrow D'
	\]
	of type III surfaces.\\
	If the MHS on $H^2(D)$ and $H^2(D')$ are split, then their period points are trivial, so the condition that $\mu$ induces an isomorphism of MHS becomes vacuous. The last claim follows.
\end{proof}
\section{Torelli for log Calabi--Yau threefolds}
\subsection{The mixed Hodge structure on $H^3(U)$}
Let $(Y, D)$ be a three-dimensional log Calabi--Yau pair with a toric model. Suppose also that $H^3(Y)=0$.
$Y$ is a (smooth) rational threefold, so that $\Pic(Y) \cong H^2(Y, \ZZ)$. 
The excision sequence in cohomology gives 
\[
\dots 	\rightarrow H_4(Y) \rightarrow H^2(D)(3) \rightarrow H_3(U) \rightarrow H_3(Y)=0.
\]
It follows that the MHS $H_3(U)$ is a quotient of $H^2(D)(3)$: more precisely, writing $K$ for the image of $H_4(Y) \cong H^2(Y) \cong \Pic(Y)$ in $H^2(D)$, we have an isomorphism
\[
(H^2(D)/ K)(3) \cong H_3(U)
\]
of mixed Hodge structures. 
The sub-structure $K \subset H^2(D)$ has pure weight $2$, so we obtain from~\eqref{eq:ext4} an extension
\begin{equation}\label{eq:ext2}
	0 \rightarrow \ZZ(3) \rightarrow H_3(U) \rightarrow (\Lambda/K)(3) \rightarrow 0
\end{equation}
which expresses $H_3(U)$ as an extension of a weight $-4$ Hodge structure by a Hodge structure of weight $-6$. Let 	
$
\bar{\psi} \colon \Lambda/K \rightarrow \CC^\times
$
be the associated extension class. We show:
\begin{lemma}\label{lem:inducedextension}
	$\bar{\psi}$ is induced by the extension class $\psi \colon \Lambda \rightarrow \CC^\times$ of \eqref{eq:ext4}.
\end{lemma}
\begin{proof}
	Clearly, $K \subset \Pic(D)$, so by Corollary~\ref{cor:periodpoint}, the restriction of $\psi \colon \Lambda \rightarrow \CC^\times$ to $K$ is trivial. It follows that the pullback of the extension~\eqref{eq:ext4} by the inclusion $K \hookrightarrow \Lambda$ is a split extension. Using that $\Hom_{\MHS}(\Lambda, \ZZ(0))=0$ (there are no morphisms preserving the Hodge filtration), the long exact sequence of Ext groups shows that the extension \eqref{eq:ext4} is the pullback of the extension \eqref{eq:ext2} via the projection $p \colon \Lambda \twoheadrightarrow \Lambda/K$. It follows that 
	\[
	{\psi}=\bar{\psi} \circ p,
	\]
	as required.
\end{proof}
In particular, this shows that the MHS $H_3(U)$, or equivalently its dual $H^3(U)$, is completely determined by $H^2(D)$. This is crucial to our argument, and fails as soon as the assumption $H^3(Y)=0$ is dropped.
\subsection{Torelli for log CY3s with a common toric model}
If a type III surface $D$ arises as the boundary of a log CY3 pair $(Y, D)$, the exact sequence~\eqref{eq:seq2} can be continued on the left:
\begin{lemma}
	Suppose that $(Y, D)$ is a log CY3 pair which admits a toric model $(Y, D) \rightarrow (\bar{Y}, \bar{D})$.
	There is an exact sequence 
	\begin{equation}\label{eq15}
		0 \rightarrow \ker \left(\Aut(Y, D) \rightarrow \Pic(Y) \right) \rightarrow \Aut^0(D^{[1]}) \xrightarrow{\psi} \bigoplus_{v \in \Sigma^{[0]}}\Hom(\Pic(D_v), \CC^\times) \rightarrow \Hom(\Lambda, \CC^\times) \rightarrow 0
	\end{equation}
\end{lemma}
\begin{proof}
	Let $N$ be the cocharacter lattice of $\bar{Y}$. Each component $\bar{D}_v$ corresponds to a vector $n_v \in N$, the first lattice point on the ray corresponding to $\bar{D}_v$.
	We have a homomorphism
	\[
	\bigoplus_{e \in \Sigma^{[1]}} \Pic(\bar{D}_e)\cong \bigoplus_{e \in \Sigma^{[1]}}\ZZ[\bar{D}_e] \rightarrow \wedge^2N
	\]
	defined by sending the basis vector $[\bar{D}_e]$ to $n_v \wedge n_w$, where $e=e(v,w)$ and the order in the wedge product is determined by the orientation of edges. It follows that 
	\begin{equation}\label{eq:gamma}
		\gamma(L_v|_{\bar{D}_e})=(L_v \cdot \bar{D}_e) n_v \wedge n_w. 
	\end{equation} 
	We claim that the sequence from \eqref{eq:seq1} can be continued to the right:
	\[
	\begin{tikzcd}
		0 \ar[r]& \Pic(\bar{D}) \ar[r] & \bigoplus_{v \in \Sigma^{[0]}} \Pic(\bar{D}_v) \ar[r, "\ell"] &\bigoplus_{e \in \Sigma^{[1]}} \Pic(\bar{D}_e) \ar[r, "\gamma"] & \wedge^2N \ar[r] &0.
	\end{tikzcd}
	\]
	Exactness at the first and second term follows immediately from Corollary~\ref{cor:periodpoint}. Since the vectors $n_v$ span $N$, exactness at the last term is also clear.
	By \cite[Proposition 3.2.7]{CoxLittleSchenck}, the cocharacter lattice of the component $\bar{D}_v$ is $N(v)=N/\langle n_v \rangle$, and the toric divisor $\bar{D}_e$ on $\bar{D}_v$ corresponding to an edge $e=e(v,w)$ corresponds to the image $\bar{n}_w$ of $n_w$ in $N(v)$. Given $L_v \in \Pic(\bar{D}_v)$, the exact sequence defining the Picard group of a toric surface~ \cite[Theorem 4.1.3]{CoxLittleSchenck} shows that 
	\[
	\sum_e (L_v \cdot \bar{D}_e) \bar{n}_w=0,
	\]
	where the sum is over all toric divisors $\bar{D}_e$ of $\bar{D}_v$, and therefore 	
	\[
	\sum_e (L_v \cdot \bar{D}_e) n_v \wedge n_w=0
	\] 
	as well. 
	Using \eqref{eq:gamma}, it now follows easily that $ \Ima(\ell) \subset \ker(\gamma)$. We obtain an induced surjective morphism $\bar{\gamma} \colon \coker(\ell) \rightarrow \wedge^2 N$ and it remains to show $\bar{\gamma}$ is an isomorphism. Since $\bar{\gamma}$ is a morphism of finitely generated $\ZZ$-modules, it is enough to show that $\coker(\ell) \cong \wedge^2 N$. By \cite[Example 3.5]{Petersen} or \cite[Corollary 8.33]{Voisin}, there is a spectral sequence 
	\[
	E_1^{p,q}=\oplus_{i=p}H^q(\bar{D}^{[i]}, \ZZ)
	\]
	converging to $H_c^{p+q}(\bar{Y} \setminus \bar{D})$.
	By \cite[Proposition 8.34]{Voisin}, the differential \[
	d_2 \colon E_1^{1,2} \rightarrow E_1^{2,2}
	\]
	is identified with $\ell$, and since $E_1^{3,2}=0$, we have $E_2^{2,2}=\coker(\ell)$. For degree reasons, all further differentials are zero, and also $E_2^{p,q}=0$ if $p+q=4$ and $(p,q) \neq (2,2)$. We therefore have an isomorphism 
	\[
	E_2^{2,2} \cong H_c^4(U, \ZZ) \cong H_2(U, \ZZ).
	\]
	Since $U=(\CC^\times)^3$, the group $H_2(U, \ZZ)$ is naturally identified with $\wedge^2 N$.\\
	For $(Y, D)$, we thus obtain an exact sequence	
	\[
	\begin{tikzcd}
		0 \ar[r]&\Lambda \ar[r]& \bigoplus_{v \in \Sigma^{[0]}} \Pic(D_v) \ar[r] &\bigoplus_{e \in \Sigma^{[1]}} \Pic(D_e) \ar[r]& N' \ar[r] &0
	\end{tikzcd}
	\]
	where $N'$ is the quotient of $\wedge^2N$ by the subgroup generated by tensors $n_v \wedge n_w$ corresponding to one-dimensional strata of $D$ which intersect the centers of the blowups in the toric model $(Y, D) \rightarrow (\bar{Y}, \bar{D})$. 
	Applying the exact functor $\Hom(-, \CC^\times)$, we obtain an exact sequence 
	\[
	0 \rightarrow \Hom(N', \CC^\times) \rightarrow \Aut^0(D^{[1]}) \xrightarrow{\phi} \bigoplus_{v \in \Sigma^{[0]}}\Hom(\Pic(D_v), \CC^\times) \rightarrow \Hom(\Lambda, \CC^\times) \rightarrow 0.
	\]
	We note that $\Hom(N', \CC^\times) \subset \Hom(\wedge^2N, \CC^\times) \cong \Hom(M, \CC^\times)=T_N$ is the subgroup of $\Aut(\bar{Y}, \bar{D})$ fixing pointwise those one-dimensional strata along which we are blowing up. It is easy to see that this group is identified with $ \ker \left(\Aut(Y, D) \rightarrow \Pic(Y) \right)$.
\end{proof}
Let $(Y, D)$ be a log CY3 pair with maximal boundary, so that $D$ a type III surface. Let $\pi \colon (Y', D') \rightarrow ({Y}, {D})$ be an interior blowup with exceptional divisor $E$ as in Lemma~\ref{lem:centersoftoricmodel}. Let $\ell$ be the minimal generator of the corresponding extremal ray of $\overline{\NE}(Y)$. The sequence 
\[
0 \rightarrow \Pic(Y) \xrightarrow{\pi^*} \Pic(Y') \xrightarrow{\cdot (-\ell)} \ZZ \rightarrow 0
\]
from \eqref{eq24} is exact on the right. Since $E \cdot \ell=1$, it has a canonical splitting given by $1 \mapsto [E]$, so we obtain a canonical isomorphism
\[
\mu \colon \Pic(Y') \cong \Pic(Y)\oplus \ZZ.
\]
If $\pi$ is the blowup of a curve $C$, the exceptional $E$ divisor meets one component of $D'$ in a section of $E \rightarrow C$, and adjacent components $D'_v$ in a (possibly empty) set of $n_v$ disjoint smooth rational curves $\{E^i_v\}$. 
\begin{definition}
	A \emph{marking} of exceptional curves is an ordering of the curves $E^i_v$, for each $v \in \Sigma^{[0]}$. 
\end{definition}
Given a marking of exceptional curves, we obtain a canonical isometry
\[
\mu_v \colon \Pic(D'_v) \cong  \Pic({D}_v)\oplus   \ZZ^{n_v}
\]
for every $v\in \Sigma^{[0]}$. \\
If $\pi$ is the blowup of a point, then the exceptional $E$ divisor meets two components of $D'$ in a smooth rational curve, and intersects all other components trivially. We obtain a canonicial isometry
\[
\mu_v \colon \Pic(D'_v) \cong  \Pic({D}_v)\oplus   \ZZ^{n_v}
\]
for every $v\in \Sigma^{[0]}$ (where $n_v$ is either $0$ or $1$).
\begin{definition}
	If $(Y, D)$ and $(Y', D')$ are two pairs with a common toric model $(\bar{Y}, \bar{D})$, obtained by blowing up curves $C_1, \dots C_r$ (resp. $C_1', \dots C_r')$, and points $p_1, \dots p_s$ (resp. $p'_1, \dots p'_s$) in a specified order, we say that the pairs are of the \emph{same combinatorial type} if 
	\[
	C_k \cdot D_e=C_k' \cdot D_e
	\]
	for all $1 \leq k \leq r$ and all one-dimensional strata $D_e$, and if the number of points $p_i$ on $D_e$ is the same is the same as the number of points $p'_i$ on $D'_e$.
\end{definition}
If $(Y, D)$ and $(Y', D')$ have the same combinatorial type, a marking of exceptional curves gives canonical isomorphisms $\mu \colon \Pic(Y) \cong \Pic(Y')$ identifying the classes of exceptional divisors and components of $D$, and canonical isometries
\[
\mu_v \colon \Pic(D_v)\cong \Pic(D'_v)
\]
preserving one-dimensional strata of $D$ and exceptional curves.
\begin{proposition}\label{prop:TorelliCommonToricModel}
	Let $(Y, D)$ and $(Y', D')$ be log CY3 pairs with a common toric model $(\bar{Y}, \bar{D})$ of the same combinatorial type. Fix markings $p_e, p_e'$ of $D$ and $D'$, and a marking of exceptional curves. For every $v \in  \Sigma^{[0]}$, let $\mu_v \colon \Pic(D_v) \rightarrow \Pic(D'_v)$ be the canonical induced isometry. Suppose that the curves $C_k$ (resp $C_k'$) blown up in the toric model are smooth rational, and that $\oplus_v \mu_v([C_k])=[C_k']$.
	\begin{enumerate}
		\item If $\phi_{(D, p_e)}=\phi_{(D', p_e')} \circ \oplus_v\mu_v$,
		then $\mu$ and $\oplus_v \mu_v$ are induced by an isomorphism 
		\[
		f \colon (Y, D, p_e) \rightarrow (Y', D', p_e').
		\]
		\item 
		If $\phi_D= \phi_{D'} \circ \oplus_v\mu_v|_{\Lambda}$, then $\mu$ and $\oplus_v \mu_v$ are induced by an isomorphism 
		\[
		f \colon (Y, D) \rightarrow (Y', D').
		\]
	\end{enumerate}
\end{proposition}
\begin{proof}
	(1) The restriction of $\oplus_v\mu_v$ to $\bigoplus_{v \in \Sigma^{[0]}} \Pic(\bar{D}_v)$ is the identity map, so tautologically, the marked period points $\phi_{(\bar{D}, p_e)}$ and $\phi_{(\bar{D}, p_e')}$ for the toric pair $(\bar{Y}, \bar{D})$ agree. By \eqref{eq15}, 
	the element $\alpha \in \Aut^0(D^1)$ sending the $p_e$ to the $p_e'$ extends to an automorphism $f_0 \colon (\bar{Y}, \bar{D}, p_e) \rightarrow (\bar{Y}, \bar{D}, p'_e)$.
	We now argue inductively. By construction $(Y, D)$ is obtained by a sequence of interior blowups 
	\[
	(Y, D)\rightarrow \dots \rightarrow (\bar{Y}, \bar{D}).
	\]
	Suppose that we have extended $f_0$ to an isomorphism $f_k \colon (Y_k, D_k, p_e) \rightarrow (Y_k', D_k', p_e')$ inducing $\mu$ and the $\mu_v$. \\
	{\bf Case 1}: The next interior blowup is the blowup of a smooth rational curve $C_{k+1}$ on $Y_k$ (respectively $C'_{k+1}$ on $Y_k'$). 
	We will show that the isomorphism $f_k$ takes the curve $C_{k+1}$ to $C_{k+1}'$.
	Let $D_0$ be the component of $D$ containing $C_{k+1}$, and fix a component $D_1$ adjacent to $D_0$. Let $D_e$ be the one-dimensional stratum $D_0 \cap D_1$.  Let $q_i$ be the (possibly empty) set of points where $C_{k+1}$ meets $D_e$. Let $E_i$ be the exceptional curves on the strict transform of $D_1$ arising from the blowup of $C_{k+1}$. 
	We use analogous notation on $D'$. 
	Since $\oplus_v \mu_v$ preserves the marked period point, we have that
	\[
	\cO_{\partial D_1}(p_e-q_i)=\phi_{D, p_e}(E_i)=\phi_{D', p'_e}\circ \oplus_v \mu_v(E_i)=\phi_{D', p'_e}(E'_i)=\cO_{\partial D'_1}(p'_e-q'_i).
	\]
	Since $f_k(p_e)=p_e'$, the isomorphism $f_k$ takes $q_i$ to $q_i'$. 
	Repeating this argument for all components of $D$ adjacent to $D_0$, we see that $f_k$ takes the points $C_{k+1} \cap  \partial D_0$ to $C'_{k+1} \cap \partial D_0'$. It follows that the curves $f_k(C_{k+1})$ and $C'_{k+1}$ meet $\partial D'_0$ in the same points. Since the curves are smooth rational, they must be equal by Lemma~\ref{lem:restriction} below.
	It follows that that $f_k$ extends to an isomorphism $f_{k+1} \colon (Y_{k+1}, D_{k+1}, p_e) \rightarrow (Y_{k+1}', D_{k+1}', p_e')$, and by construction, this extension induces $\mu$ and $\mu_v$.\\
	{\bf Case 2}: The next interior blowup is the blowup of a point $q$ (respectively $q'$).
	We need to show that $f_k(q)=q'$. 
	The exceptional divisor meets two components $D_0$ and $D_1$ in a smooth rational curve $E_0$ (resp. $E_1$). Using again that $\oplus_v \mu_v$  preserves the period point 
	we have
	\[
	\cO_{\partial D_0}(p_e-q)=\phi_{D, p_e}(E_0)=\phi_{D', p'_e}\circ \oplus_v \mu_v(E_0)=\phi_{D', p'_e}(E'_0)=\cO_{\partial D'_0}(p'_e-q').
	\]
	Since $f_k(p_e)=p_e'$, the isomorphism $f_k$ takes $q$ to $q'$, as required.\\
	(2) In the unmarked case, choose arbitrary markings $p_e$ and $p_e'$ and let ${\phi}_{D, p_e}$ and ${\phi}_{D', p_e'}$ be the associated marked period points. By the exact sequence~\eqref{eq15}, we may change the marking on $D$ to a different set of points $r_e$ such that the marked period points satisfy ${\phi}_{D', p_e'} \circ  \oplus_v\mu_v={\phi}_{D, r_e}$. We have therefore reduced to the marked case.
\end{proof}
\begin{lemma}\label{lem:restriction}
	Let $(Y, D)$ be a Looijenga pair and $C$ be a smooth rational curve on $Y$. 
	Then the restriction map 
	\[
	H^0(Y, \cO_Y(C)) \rightarrow H^0(D, \cO_D(C))
	\]
	is injective. In other words, $C$ is uniquely determined within its linear equivalence class by its restriction to $D$.
\end{lemma}
\begin{proof}
	Twisting the standard exact sequence of the divisor $D$ by $\cO(C)$ and taking cohomology we see that it is enough to show that $h^0(\cO_Y(C-D))=0$, i.e that $A=C-D$ is not effective. Adjunction for $C \cong \PP^1$ gives $C \cdot (C-D)=-2$, so if $A$ is effective, then $C$ is fixed in the linear system $A$. This implies that $-D$ is effective, a contradiction. 
\end{proof}
\subsection{Generic Torelli}
We now work towards proving Theorem~\ref{thm:GenericTorelli}. Our strategy will be to reduce to the case where the two pairs $(Y, D)$ and $(Y', D')$ have a common toric model.
We start with two lemmas. The first gives a criterion for recognizing when a given pair $(Y, D)$ is toric. The second gives a criterion for when two toric pairs are isomorphic.  
\begin{lemma}\label{lem:ToricChararcterization}
	Let $(Y, D)$, $D=\sum_i D_i$ be a log CY pair and $(\bar{Y}, \bar{D})$, $\bar{D}=\sum_i \bar{D}_i$ be a toric log CY pair. Suppose that there is an isomorphism $\mu \colon \Pic(\bar{Y}) \rightarrow \Pic({Y})$ such that $\mu(\bar{D}_i)={D}_i$ for all $i$. Then $(Y, D)$ is in fact a toric pair.
\end{lemma}
We deduce this from \cite[Theorem 1.2]{BrownMcKernanSvaldiZong}.
We first recall the notion of complexity of a log pair $(X, \Delta)$:
Let $X$ be a variety of dimension $n$, and let $(X, \Delta)$ be a log pair. A decomposition of $\Delta$ is an expression $\Delta=\sum_{i=1}^k a_iD_i$, where $a_i \geq 0$ and the $D_i$ are Cartier divisors. The complexity of the decomposition is $c=n+r-d$, where $d=\sum a_i$ and $r$ is the rank of the vector space spanned by the $D_i$ in $\Pic(X)$. 
The complexity $c=c(X, \Delta)$ is the infimum of the complexity of any decomposition of $\Delta$. 
Well-known results in toric geometry show that if $(X, \Delta)$ is a toric variety and $\Delta$ is its toric boundary, then $d=n+r$, so $c=0$ for any decomposition.
Consider the log Calabi-Yau pair $(Y, D)$, where $D=\sum_i D_i=\sum_i \mu(\bar{D}_i)$. Since $\Pic(\bar{Y}) \rightarrow \Pic(Y)$ is an isomorphism, we have $c(Y, D)=c(\bar{Y}, \bar{D})=0$. \cite[Theorem 1.2]{BrownMcKernanSvaldiZong} then implies that there is a divisor $D'$ such that $(Y, D')$ is a toric pair, with $D' \geq D$ and $D'$ supported on $D$.
Since both $D'$ and $D$ are anticanonical, it follows that $D=D'$, as required.
\begin{lemma}\label{lem:TVdetermined}
	Let $(\bar{Y}, \bar{D})$ and $(\bar{Y'}, \bar{D'})$ be toric pairs with $\Sigma_D \cong \Sigma_D'$ as simplicial complexes
	Suppose that there are isomorphisms
	\[
	\mu \colon \Pic(\bar{Y}) \rightarrow \Pic(\bar{Y}'), \qquad \mu_v \colon \Pic(\bar{D}_v) \rightarrow \Pic(\bar{D}_v'), \; v \in \Sigma^{[0]}
	\]
	which are compatible, in the sense that the diagram 
	\[
	\begin{tikzcd}\label{diagram2}
		\Pic(\bar{Y}) \ar[r, "\mu"] \ar[d]& \Pic(\bar{Y}') \ar[d]\\
		\Pic(\bar{D}_v) \ar[r, " \mu_v"]&\Pic(\bar{D}_v')&
	\end{tikzcd}
	\]
	commutes for every $v \in \Sigma^{[0]}$. Suppose moreover that 
	\begin{itemize}
		\item $\mu(\bar{D}_v)=\bar{D}'_{v}$ and $\mu_v(\bar{D}_e)=\bar{D}'_{e}$.
		\item $\mu$ and the $\mu_v$ preserve the intersection product.
	\end{itemize}
	Then $\mu, \mu_v$ are induced by an isomorphism $f \colon (\bar{Y}, \bar{D}) \rightarrow (\bar{Y'}, \bar{D'})$.
\end{lemma}
\begin{proof}
	Let $\Sigma_{\bar{Y}}$ be the fan of $\bar{Y}$, with underlying lattice $N$. For every component $D_i$ of $D$, denote $v_i$ the minimal generator of the corresponding ray in $\Sigma$. Choose $3$ components $D_1, D_2, D_3$ meeting in a point. Since $\bar{Y}$ is smooth, there is a unique isomorphism $\phi \colon N \rightarrow \ZZ^3$ sending $v_1, v_2, v_3$ to the standard basis of $\ZZ^3$, and we define $\Sigma=\phi(\Sigma_{\bar{Y}})$, so that $\bar{Y} \cong Y_{\Sigma}$. Consider now a component $D_4$ of $D$ which is adjacent to two of $\{D_1, D_2, D_3\}$, say $D_1$ and $D_2$, and let $C=D_1 \cap D_2$.
	Standard results on intersection products on toric varieties \cite[Proposition 6.4.5]{CoxLittleSchenck} show that we have an equality 
	\begin{equation}\label{eq:7}
		v_3+c_1v_1+c_2v_2+v_{4}=0
	\end{equation}
	where $c_1=C \cdot D_1=(C^2)_{D_2}$ and $c_2=C \cdot D_2=(C^2)_{D_1}$. 
	Similarly, let $\Sigma_{\bar{Y}'}$ be the fan of $\bar{Y}'$, define an isomorphism $\phi'$ in a similar fashion, and let $\Sigma'$ be the corresponding fan in $\ZZ^3$.
	Since $\mu$ preserves the dual complex and the $\mu_v$ are isometries preserving one-dimensional strata of $D$ we have $\phi'(v'_4)=\phi(v_4)$. 
	Repeating, we find that $\phi'(v'_i)=\phi(v_i)$ for all $i$, so the fans $\Sigma$ and $\Sigma'$ have the same rays. 
	The rays $v_i$ and $v_j$ span a cone of $\Sigma_{\bar{Y}}$ if and only if $D_i$ and $D_j$ meet along a one-dimensional stratum. Using the commutative diagram in Lemma~\ref{lem:TVdetermined}, we see that this happens if and only $D_i'$ and $D_j'$ meet along a one-dimensional stratum, i.e $v_i'$ and $v_j'$ span a cone in $\Sigma_{\bar{Y}'}$. 
	Finally, the rays $v_i, v_j, v_k$ span a cone if and only if $D_i \cdot D_j \cdot D_k=1$. Since $\mu$ preserves the intersection product, this happens if and only if $v_i', v_j', v_k'$ span a cone. This implies that the fans $\Sigma$ and $\Sigma'$ have the same cones, and therefore $Y_{\Sigma}=Y_{\Sigma'}$. The result follows.
\end{proof}
Recall that on a smooth projective surface $D$, the intersection product defines a quadratic form on $\Pic(D)$. 
On a smooth projective threefold $Y$, the intersection product defines a cubic form on $\Pic(Y)$.
We are now ready to state the global Torelli theorem. 
\begin{theorem}\label{thm:GlobalTorelliv1}
	Let $(Y, D)$ and $(Y', D')$ be log CY 3 pairs with maximal boundary such that $\Sigma_D \cong \Sigma_{D'}$ as CW complexes, and $H^3(Y)=H^3(Y')=0$. Suppose that $(Y, D)$ admits a toric model 
	\[
	(Y, D) \rightarrow (Y^{n}, D^{n})  \rightarrow \dots \rightarrow (Y^0, D^0)=(\bar{Y}, \bar{D}).
	\]
	Suppose further that there are isomorphisms
	\[
	\mu \colon \Pic(Y) \rightarrow \Pic(Y'), \qquad \mu_v \colon \Pic(D_v) \rightarrow \Pic(D_v'), \; v \in \Sigma^{[0]}
	\]
	which are compatible, in the sense that the diagram 
	\begin{equation}\label{diagram}
		\begin{tikzcd}
			\Pic(Y) \ar[r, "\mu"] \ar[d]& \Pic(Y') \ar[d]\\
			\Pic(D_v) \ar[r, "\mu_v"]&\Pic(D'_v)&
		\end{tikzcd}
	\end{equation}
	commutes for every $v \in \Sigma^{[0]}$, and
	\begin{enumerate}
		\item $\mu(D_v)=D'_{v}$ and $\mu_v(D_e)=D'_{e}$.
		\item $\mu_v$ preserves the quadratic form on $\Pic(D_v)$ defined by the intersection product.
		\item $\mu$ preserves the cubic form on $\Pic(Y)$ defined by the intersection product.
		\item $\mu(\Nef(Y))=\Nef(Y')$
		\item $\mu(\Eff(Y))=\Eff(Y')$
		\item $\oplus_v \mu_v$ preserves the period point i.e. $\phi_{D'} \circ \oplus_v \mu_v|_{\Lambda}=\phi_D$
	\end{enumerate}
	Then there is an isomorphism of pairs $f \colon (Y, D) \rightarrow (Y', D')$ inducing $\mu$ and the $\mu_v$.
\end{theorem}
\begin{proof}
	We will use Mori's classification of threefold extremal rays to construct a toric model for $(Y', D')$. Let $E$ be the exceptional divisor of $\pi \colon Y \rightarrow Y^n$. Since $H^3(Y)=0$, the image $\pi(E)$ is a smooth rational curve $C_n$ contained in a component $D^n_{v_0}$ of $D^n$. 
	Denote
	\begin{itemize}
		\item  $R$ the $K$-negative extremal ray of $\overline{\NE}(Y)$ contracted by $\pi \colon Y \rightarrow Y^n$.
		\item $\ell \in \overline{\NE}(Y)$ the primitive generator of $R$
		\item  $S$ the extremal ray of $\Eff(Y)$ spanned by $E$.
		\item  $H$ a nef divisor giving the contraction $\pi$.
	\end{itemize}
	Note that $S\cdot R<0$ by negativity of contraction.
	Since $\pi$ is divisorial, $H$ is big, and therefore $H^3>0$.
	Since $\mu$ preserves the intersection product and the nef cone, $\mu(H)$ is also big and nef, and some multiple $H'$ defines a birational map $\pi' \colon Y' \rightarrow Y'^{n}$. 
	By duality, and using that $\mu$ preserves the canonical class, $\pi'$ is the contraction of the $K$-negative extremal ray $\mu(R)$ of $\overline{\NE}(Y')$. We set
	\begin{itemize}
		\item  $R'=\mu(R)$
		\item $\ell' \in \overline{\NE}(Y)$ the primitive generator of $R'$.
		\item  $S'$ the extremal ray of $\Eff(Y')$ spanned by the exceptional divisor $E'$ of $\pi'$. 
	\end{itemize}
	We claim that $S'=\mu(S)$. By negativity of contraction, we have $S' \cdot R' < 0$, and by duality, $\mu(S) \cdot R' <0$ as well. Since $S'$ and $\mu(S)$ are extremal rays of $\Eff(Y')$, they both contain the class of an irreducible divisor. 
	However, the only irreducible divisor negative on $R'$ is $E'$, so that $\mu(S)$ and $S'$ contain $[E']$ and hence $\mu(S)=S'$. In particular, we have that 
	\begin{equation}\label{eq23}
		\mu([E])=r[E'] 
	\end{equation}
	for some $r \in \QQ_{>0}$. By Lemma~\ref{lem:primitive}, both $[E]$ and $[E']$ are primitive generators of their extremal rays, so $r=1$.\\
	{\bf Case 1:} Suppose first that $\pi$ is a blowup of a smooth curve on a component $D^n_{v_0}$ of $D^n$ meeting one-dimensional strata transversely. For $v \neq v_0$, the intersection $D_v \cdot E$ is a (possibly empty) set of disjoint $(-1)$-curves $\{E^i_v\}$ on $D_v$. 
	For $v = v_0$, the intersection $D_v \cdot E$ is a section $C$ of $\pi \colon E \rightarrow C_n$. Denoting $i \colon H_2(D_v) \rightarrow H_2(Y)$, we see that $i([C]) \notin R$ , but $i([E^i_v]) \in R$.
	We claim that $\pi' \colon Y' \rightarrow Y'^n$ is also the blowup of a smooth curve (Type 1 in Mori's classification). Indeed, by \cite[Lemma 3.20, 3.21]{Mori}, $\pi$ is of type 1 if and only if $H \cdot E \neq 0 \in H_2(Y)$. 
	By assumption, $\mu$ preserves the intersection product, so we have $\mu(H) \cdot \mu(E) \neq 0 \in H_2(Y')$.
	Since $[H']$ and $[E']$ are positive multiples of $\mu([H])$ and $\mu([E])$, we obtain that $H' \cdot E' \neq 0 \in H_2(Y')$ as well, 
	so that $\pi'$ is the blowup of a smooth curve $C_n'$, which is smooth rational since $H^3(Y')=0$.\\
	The exceptional divisor $E'$ is not a component of $D'$ (as the isomorphism $\mu$ preserves the boundary components of $D$), and similarly, $\ell'$ is not numerically equivalent to a one-dimensional stratum of $D'$. 
	Since $\ell'$ sweeps out a divisor, we have $D'_v \cdot \ell' \geq 0$ for every component $D'_v$ of $D'$. On the other hand, $D' \cdot \ell'=1$ by duality, so $E'$ meets exactly one component $D'_{v'_0}$ of $D'$ in a section $C'$ of $\pi' \colon E' \rightarrow C'_n$, and every other component $D'_v$ in a union $\{{E'}^i_v\}$ of fibers of $E'$. The ${E'}^i_v$ are $(-1)$-curves on $D'_v$, and the restriction of $\pi'$ to $D'_v$ contracts these $(-1)$-curves. It follows that $D'^n$ is again a normal crossing divisor with smooth components. Moreover, $C_n'$ meets the boundary of $D'_{v'_0}$ transversely.
	Since $\pi'$ is birational and contracts a ruled surface to a smooth curve, we have 
	\[
	-K_{Y'^{n}} = \pi'_*\pi'^*(-K_{Y'^{n}})=\pi'_*(-K_{Y'}+E')=\pi'_*(D'),
	\]
	so that $D'^n=\pi'(D') \in |-K_{Y'^{n}}|$, and therefore, $(Y'^n, D'^n)$ is log CY3, with $D'^n$ a type III surface .\\
	We wish to show that the assumptions of the theorem are satisfied for the pairs $(Y^n, D^n)$ and $(Y'^n, D'^n)$. 
	Commutativity of \eqref{diagram} shows that 
	\[
	\mu_v([E|_{D_v}])=E'|_{D_v}
	\]
	for all $v \in \Sigma^{[0]}$. If $v=v_0$, then $[E|_{D_{v_0}}]=[C]$, and in particular, $i_*([C]) \notin R$.
	Using the dual diagram to~\eqref{diagram}:
	\[
	\begin{tikzcd}
		H_2(Y) \ar[r, "(\mu^*)^{-1}"] & H_2(Y') \\
		H_2(D_v) \ar[u,"i_*"]\ar[r, " \mu_v"]&H_2(D'_v) \ar[u, "i_*"]
	\end{tikzcd}
	\]
	we see that $i_*([E'|_{D'_{v_0}}]) \notin R'$. This is only possible if $D'_{v_0}$ is the component of $D'$ meeting $E$ in a section, i.e. $v_0=v_0'$ and 
	\begin{equation}\label{eq18}
		\mu_{v_0}([C])=[C']
	\end{equation}
	For the other components, we have 
	\[
	[E|_{D_v}]=\sum_i [E^i_v], \qquad [E'|_{D'_v}]=\sum_i [{E'}_v^i]
	\]
	and we claim that $\mu_v$ maps the set $\{[E^i_v]\}$ to the set $\{[{E'}_v^i]\}$.
	Indeed, since $\mu$ preserves the nef cone, $\mu_v$ sends some ample class to an ample class. Since $\mu_v$ is an isometry, it then follows from Riemann-Roch that $\mu_v([E^i_v]) \in \Pic(D_v')$ is represented by an effective divisor of self-intersection $(-1)$. The image of this class under $i_*$ is in the extremal ray $R'$, so we see that 
	$\mu([E^i_v])$ is equal to $[{E'}^j_v]$ for some $j$. 
	Up to changing the  marking (ordering) of the ${E'}^i_v$, we therefore have 
	\begin{equation}\label{eq67}
		\mu_v([E^i_v])=[{E'}^i_v]
	\end{equation}
	for all $i$. By \eqref{eq24}, we have an exact sequence
	\begin{equation}\label{eq224}
		0 \rightarrow \Pic(Y^n) \xrightarrow{\pi^*} \Pic(Y) \xrightarrow{ \cdot (-\ell)} \ZZ \rightarrow 0
	\end{equation}
	where exactness on the right follows from $E \cdot \ell=-1$, and the minus sign was inserted for convenience. In particular, the image of $\pi^*$ is identified with $\ell^\perp$. We have a similar exact sequence relating $\Pic(Y')$ and $\Pic(Y'^{n})$.
	Define the isomorphism
	\[
	\mu^n \colon \Pic(Y^n) \rightarrow \Pic(Y'^{n}), \qquad \mu^n=(\pi'^*)^{-1} \circ \mu \circ \pi^*
	\]
	Using Lemma~\ref{lem:conespullback}, we see that $\mu^n$ preserves the effective cone and the nef cone. Since $\pi^*$ and $\mu$ preserves the intersection product and the classes of the components of $D$, the same holds for $\mu^n$. 
	The inclusion $i \colon 1 \mapsto [E]$ defines a splitting of \eqref{eq24}, so that 
	\[
	\Pic(Y^n)\oplus \ZZ \xrightarrow{(\pi^*, i)} \Pic(Y)
	\]
	an isomorphism. Since $\mu([E])=[E']$, we see that under this identification (and a similar one on $Y'$), the isomorphism $\mu^n \oplus 1$ becomes $\mu$.\\
	Turning to the boundary, the induced contraction $\pi \colon D_v \rightarrow D^n_v$ is the contraction of a $K$-negative extremal face, which is the contraction of $k_v$ disjoint $(-1)$-curves $\{E^i_v\}$ on $D_v$. In what follows we write $E_i$ for $E_v^i$. Our ordering of the exceptional curves gives an exact sequence 
	\begin{equation}\label{eq26}
		0 \rightarrow \Pic(D^n_v) \xrightarrow{\pi^*} \Pic(D_v) \xrightarrow{\rho} \ZZ^{k_v} \rightarrow 0.
	\end{equation}
	where 
	\[
	\rho(D)=( -D\cdot E_1,..., -D \cdot E_{k_v})
	\]
	
	We have a similar exact sequence for $Y'$.
	By \eqref{eq67}, $\mu_v$ induces an isomorphism $\langle E_1, \dots E_{n_v}\rangle^\perp \rightarrow \langle E_1, \dots E_{n_v}\rangle^\perp$.
	For every $v \in \Sigma^{[0]}$, we can therefore define the isometry
	\[
	\mu^n_v \colon \Pic(D^n_v) \rightarrow \Pic(D'^n_v), \qquad \mu^n_v =({\pi'}^*)^{-1} \circ \mu_v \circ \pi^*
	\]
	and by \eqref{eq18}, we have 
	\[
	\mu^n_v([C_n])=[C_n']
	\]
	As before, the splitting of \eqref{eq26} given by $j_i \colon e_i \mapsto E_i$ induces an isomorphism 
	\[
	\Pic(D^n_v) \oplus \ZZ^{k_v} \xrightarrow{(\pi^*, \oplus_i j_i)}\Pic(D_v)
	\]
	and under this identification (and a similar one on $Y'$), the isomorphism $\mu_v^n \oplus {\bf 1}$ becomes $\mu_v$.\\
	Since $\mu^n$ and $\mu_v^n$ are induced by $\mu$ and $\mu_v$, the isometry $\oplus_v \mu_v^n$ preserves the period point, and we obtain a commutative diagram~\ref{diagram} for the pairs $(Y^n, D^n)$ and $(Y'^{n}, D'^{n})$. Since $\pi^*$ and $\mu_v$ preserve the classes of one-dimensional strata, the same holds for $\mu_v^n$. \\
	{\bf Case 2:} Suppose now that $\pi$ is the blowup of a point on a one-dimensional stratum of $D^n$. We keep the notation from before. In this case, there are two components $D_{v_0}$ and $D_{v_1}$ of $D$ such that $E_0=[E|_{D_0}]$ and $E_1=[E|_{D_1}]$ are the classes of $(-1)$-curves, and $[E|_{D_v}]=0$ for all other components. We wish to show that $\pi'$ is also the blowup of a smooth point. Since $\mu$ preserves the intersection product, we have $E' \cdot H'=0$, which rules out type 1 on Mori's list, and by computing $E \cdot K_Y^2=4$, we can use the table in Figure~\ref{table} to rule out all other cases except for 2). As before, $E'$ is not a component of $D'$, and $\ell'$ is not linearly equivalent to a one-dimensional stratum. By commutativity of the diagram, the classes $[E'|_{D'_{v_0}}] \in \Pic(D'_0)$ and $[E'|_{D'_{v_1}}] \in \Pic(D'_{v_1})$ are represented by an effective divisor of self-intersection $(-1)$. Since both have image $\ell$ under $i_*$, they must in fact be represented by $(-1)$-curves, and $\pi'$ blows down $E'$ to a smooth point on the one-dimensional stratum $D'_{v_0} \cap D'_{v_1}$.
	The rest of the argument is now analogous.\\
	We have shown that the assumptions of the theorem are satisfied for the log Calabi--Yau pairs $(Y^n, D^n)$ and $({Y'}^n, {D'}^n)$. 
	By induction, there is a composite 
	\[
	(Y', D') \rightarrow (Y'^{n}, D'^{n}) \rightarrow \dots \rightarrow (Y'_0, D'_0)=(\bar{Y'}, \bar{D'})
	\]
	where each map is an interior blowup, and the pairs $(\bar{Y}, \bar{D})$ and $(\bar{Y'}, \bar{D'})$ satisfy the assumptions of the Theorem. Since $(\bar{Y}, \bar{D})$ is a toric pair, 
	Lemma~\ref{lem:ToricChararcterization} shows that $(\bar{Y'}, \bar{D'})$ is also a toric pair, and Lemma~\ref{lem:TVdetermined} shows that $\mu^0$ and $\mu_v^0$ are induced by an isomorphism $(\bar{Y}, \bar{D}) \cong (\bar{Y}', \bar{D}')$. It follows that we may assume that
	$(Y, D)$ and $(Y', D')$ have a common toric model $(\bar{Y}, \bar{D})$ of the same combinatorial type. By construction $\mu_v$ is the canonical isometry associated to our choice of marking of exceptional curves, we have $\oplus_v \mu_v([C_k])=[C_k']$ and $\oplus_v \mu_v$ preserves the period point.
	We conclude using Proposition~\ref{prop:TorelliCommonToricModel}.
\end{proof}
The long list of assumptions in Theorem~\ref{thm:GlobalTorelliv1} on the isomorphisms $\mu$ and $\mu_v$ may seem daunting, but the crucial insight is that almost all these assumptions are satisfied if $\mu$ and $\mu_v$ are induced by parallel transport. We first use a standard Hilbert scheme argument to show that parallel transport preserve the generic nef and generic effective cone.
\begin{lemma}\label{lem:coneparalleltransport}
	Let $\pi \colon Y \rightarrow B$ be a smooth family of projective varieties over a connected base $B$ with $R^i\pi_*\cO_Y=0$ for $i >0$. Let $b_1, b_2 \in B$ be very general points, and let $\gamma \colon [0,1] \rightarrow B$ be a path (in the analytic topology) from $b_1$ to $b_2$.
	Then parallel transport $\gamma_* \colon H^2(Y_{b_1}, \RR) \rightarrow H^2(Y_{b_2}, \RR)$ satisfies
	\begin{align*}
		\gamma_*(\Nef(Y_{b_1}))&=\Nef(Y_{b_2})\\
		\gamma_*(\Eff(Y_{b_1}))&=\Eff(Y_{b_2}).
	\end{align*}
\end{lemma}
\begin{proof}
	The vanishing of $R^i\pi_*\cO_Y$ and cohomology and base change implies that $H^i(\cO_{Y_b})=0$ for $i>0$ and $b\in B$, and therefore $N^1(Y_b)_\RR \cong H^2(Y_b, \RR)$ for all $b \in B$. 
	We first claim that the family $\Nef(Y_b)$ is locally constant away from a countable union of closed proper subvarieties of $B$. To see this, we argue similarly to \cite[Proposition~12.2.5]{KollarMoriFlips}. We consider the relative Hilbert scheme of $\pi \colon Y \rightarrow B$. It has countably many components and each component is projective over $B$. Let $Z_i$ be the images of those components that don't dominate $B$. Note that each $Z_i$ is a closed proper subvariety of $B$. Let now $b \in B^\circ:=B \setminus \bigcup_i Z_i$ and consider a contractible analytic neighbourhood $b \in U \subset B^\circ$ over which the local system $R^2\pi_*\ZZ$ becomes trivial, so we can identify
	\begin{equation}\label{eq12}
		H^2(Y_t, \ZZ) \cong H^2(Y_b, \ZZ)
	\end{equation}
	for all $t \in U$. This identification allows us to view the nef cones $\Nef(Y_t)$ as a family of cones in the fixed vector space $H^2(Y_b, \RR)$.\\
	Consider an effective divisor $D \subset Y_b$. By construction, the component of the Hilbert scheme containing $D$ dominates $B$, so possibly after shrinking $U$, we have a codimension $1$ subscheme $\cD \rightarrow U$, flat over $U$ extending $D$, and $[\cD_t]=[D]$ under the identification $H^2(Y_t, \ZZ) \cong H^2(Y_b, \ZZ)$ given by \eqref{eq12}. It follows that the classes in $H^2(Y_t, \ZZ)$ represented by effective divisors are the same for all $t \in U$, so that the cones $\Eff(Y_t)$ are identified for all $t \in U$. The same argument works for a curve $C \subset Y_b$, so that $\NE(Y_t)$, and dually $\Nef(Y_t)$ are identified for $t \in U$. In particular, parallel transport along a path contained entirely in $U$ preserves the nef cone and the effective cone.\\
	Consider now the given path $\gamma \subset B$. By assumption, $b_1, b_2$ are very general, so contained in $B^\circ$. Since $B \setminus B^\circ$ is a countable union of algebraic subvarieties, it has real codimension at least two, so we can homotope $\gamma$ to a path $\gamma'$ entirely contained in $B^\circ$. The isomorphism on cohomology induced by parallel transport is independent of the homotopy class, so the result follows by covering $B^\circ$ with open sets $U$ as above. 
\end{proof}
\begin{proposition}\label{prop:genericpairs}
	Let $\pi \colon (\cY, \cD) \rightarrow B$ be a family of log Calabi--Yau threefolds with maximal boundary and suppose that 
	$R^i\pi_*\cO_\cY=0$
	for $i>0$. Let $\gamma \colon [0,1] \rightarrow B$ be a path from $b_1$ to $b_2$ (in the analytic topology). Then parallel transport along $\gamma$ induces an isomorphism $\Pic(Y_{b_1}) \cong \Pic(Y_{b_2})$ and compatible isomorphisms $\Pic(D_v) \cong \Pic(D_v')$ such that conditions $(1), (2)$ and $(3)$ in Theorem~\ref{thm:GlobalTorelliv1} are satisfied. \\
	Suppose further that the points $b_1, b_2$ are very general. Then $(4)$ and $(5)$ are also satisfied.
\end{proposition}
\begin{proof}
	We write $(Y_i, D_i)$ for $(\cY_{b_i}, \cD_{b_i})$ to declutter notation.
	Our definition of a family of pairs requires that the anticanonical divisor has the same combinatorial type for all fibers. In particular, parallel transport induces an isomorphism between the dual complexes $\Sigma_{D_{1}}$ and $\Sigma_{D_2}$ as CW complexes.
	After passing to an etale cover of $B$, we may assume that monodromy acts trivially on the dual complex, so that for each $v \in \Sigma^{[0]}$, we have induced families $\pi|_{\cD_v} \colon \cD_v \rightarrow B$. 
	Parallel transport in the local systems $R^2\pi_*\ZZ$ and $R^2(\pi|_{\cD_v})_*\ZZ$ gives a commutative diagram 
	\begin{equation}\label{diagram4}
		\begin{tikzcd}
			H^2(Y_{1}) \ar[r, "\mu"] \ar[d]& H^2(Y_{2}) \ar[d]\\
			H^2(({D_1})_{v}) \ar[r, "\mu_v"]&H^2(({D_2})_v)&
		\end{tikzcd}
	\end{equation}
	as in Theorem~\ref{thm:GlobalTorelliv1}, such that $\mu$ maps components of $D_1$ to the corresponding components of $D_2$, and $\mu_v$ maps one-dimensional strata of $D_1$ to the corresponding one-dimensional strata of $D_2$. 
	By cohomology and base change, the vanishing of $R^i$ shows that $H^i(\cO_{Y_i})=0$ for $i>0$, so that $H^2(Y_i, \ZZ) \cong \Pic(Y_i)$. The components of $D_i$ are rational surfaces by Proposition~\ref{lem:logcy3}, so similarly, $H^2(({D_i})_{v}) \cong \Pic({D_i})_{v})$ and we can view $\mu_v$ as an isomorphism between the Picard groups. It follows that $(1)$ holds.
	Parallel transport preserves the cup product on cohomology, so that $(2)$ and $(3)$ hold. \\
	If $b$ and $b'$ are very general points in $B$, Lemma~\ref{lem:coneparalleltransport} shows that $(4)$ and $(5)$ are satisfied as well. 
\end{proof}
We can now deduce our main theorem.
\begin{corollary}\label{cor:main}
	Let $(\cY, \cD) \rightarrow B$ be a family of log Calabi--Yau threefolds, and write $\mathcal{U}=\cY \setminus \cD$. Let $b_1, b_2 \in B$ and suppose that parallel transport along some path $\gamma$ from $b_1$ to $b_2$ induces an isomorphism $\mu_U \colon H^3(\mathcal{U}_{b_1}, \ZZ) \cong H^3(\mathcal{U}_{b_2}, \ZZ)$ of mixed Hodge structures. Suppose further that
	\begin{enumerate}
		\item $(\cY_{b_1}, \cD_{b_1})$ admits a toric model, and every $1$-dimensional center in the toric model is a smooth rational curve. 
		\item The points $b_1, b_2 \in B$ are very general.
	\end{enumerate}
	Then there is an isomorphism of pairs $f \colon (\cY_{b_1}, \cD_{b_1}) \cong (\cY_{b_2}, \cD_{b_2})$ with $(f|_{\mathcal{U}_{b_1}})_*=\mu_U$.
\end{corollary}
\begin{proof}
	We write $(Y_i, D_i)$ for $(\cY_{b_i}, \cD_{b_i})$ to declutter notation. Our definition of a family of pairs requires that the anticanonical divisor has the same combinatorial type for all fibers. In particular, parallel transport induces an isomorphism between the dual complexes $\Sigma_{D_{1}}$ and $\Sigma_{D_2}$. Since $Y_{1}$ is obtained by successively blowing up a smooth toric variety in smooth rational curves, we have $H^3(Y_{1})=0$, and therefore $H^3(Y_{2})=0$ as well. 
	By Proposition~\ref{prop:genericpairs}, parallel transport along $\gamma$ induces compatible isomorphisms $\mu \colon \Pic(Y_1) \cong \Pic(Y_2)$ and $\mu_v \colon \Pic((D_1)_v) \cong \Pic((D_2)_v)$ satisfying the conditions $(1), (2)$ and $(3)$ of Theorem~\ref{thm:GlobalTorelliv1}. Since $b_1, b_2$ are very general, $(4)$ and $(5)$ are satisfied as well.
	It remains to check that condition $(6)$ is satisfied. \\
	Recall from \eqref{eq:ext2} that the MHS $H_3(U)$ is an extension of the weight $-4$ Hodge structure $(\Lambda/K)(3)$ by the weight $-6$ Hodge structure $\ZZ(3)$. 
	The $\mu_v$ give rise to an isomorphism 
	\[
	\oplus_v\mu_v|_{\Lambda} \colon \Lambda/K \rightarrow \Lambda'/K'.
	\] 
	This isomorphism agrees with the isomorphism on weight $4$ pieces induced by $\mu_U$. Since $\mu_U$ is an isomorphism of MHS, we have
	\[
	\bar{\psi}=\bar{\psi}' \circ \oplus_v\mu_v|_{\Lambda}
	\]
	where $\bar{\psi}$ and $\bar{\psi'}$ denote the extension classes of $H_3(U_1)$ and $H_3(U_2)$. 
	By Lemma~\ref{lem:inducedextension}, $\bar{\psi}$ and $\bar{\psi}'$ are identified with the period points $\phi_D$ and $\phi_{D'}$.
	Theorem~\ref{thm:GlobalTorelliv1} now gives an isomorphism of pairs $f \colon (Y_1, D_1) \rightarrow (Y_2, D_2)$ inducing $\mu$ and the $\mu_v$, and therefore also $\mu_U$, as required.
\end{proof}
\subsection{Mirror Symmetry}\label{sec:mirrorsymmetry}
We show that families of log Calabi--Yau pairs covered by our theorem arise in mirror symmetry in the context of the Fano--LG correspondence (see \cite{MirrorSymmetryFano} for an introduction). We will follow the setup used in~\cite{Przyjalkowski, Doran}.
Let $X$ be a Fano manifold of dimension $n$, let $L \in \Pic(X)$, and let $f \in \CC[x_1^{\pm}, \dots, x_n^{\pm}]$ be a Laurent polynomial in $n$ variables. We define its \emph{classical period} $\pi_f(t) \in \CC[[t]]$:
\begin{align*}
	\pi_f(t)&=\left( \frac{1}{2\pi i} \right)^n\int_{|x_1|=\dots=|x_n|=1}\frac{1}{1-tf}\frac{dx_1}{x_1} \dots \frac{dx_n}{x_n}\\
	&=\sum_{k \geq 0} c_0(f^k)t^k
\end{align*}
where we have used the Cauchy residue theorem $n$ times in the second equality, and $c_0$ denotes the constant term of a Laurent polynomial. 
Following~\cite{Doran}, we define
\begin{definition}\label{def:mirror}
	Let $X$ be a Fano manifold and $L \in \Pic(X)$. A Laurent polynomial $f$ is a \emph{toric Landau--Ginzburg model} for $(X, L)$ if the following conditions are satisfied:
	\begin{enumerate}
		\item (Period condition) $\pi_f(t)$ is equal to $\tilde{I}^{(X, L)}(t)$, where $\tilde{I}^{(X, L)}(t)$ is a generating function for genus zero Gromov--Witten invariants of $X$ (see \cite[p.5]{Przyjalkowski} for a definition of $\tilde{I}^{(X, L)}(t)$).
		\item (Calabi--Yau condition): Viewing $f$ as a function $f \colon (\CC^\times)^n \rightarrow \CC$, there is a relative compactification $\tilde{f} \colon U_f \rightarrow \CC$ whose total space $U_f$ is a (noncompact) Calabi--Yau variety. 
		\item (Toric condition): Let $P=\Newt(f)$ be the Newton polytope of $f$, and let $X_P$ be the toric variety defined by the spanning fan of $P$\footnote{This means that $X_P$ has moment polytope $P^\vee$, the polar dual of $P$.}. Then there is a degeneration of $X$ to $X_P$.
	\end{enumerate}
	If $f \colon(\CC^\times)^n \rightarrow \CC$ extends to a morphism $\tilde{f} \colon (Y_f, D_f) \rightarrow (\PP^1, \infty)$ with $(Y_f, D_f)$ a log Calabi--Yau pair\footnote{we recall our standing assumption that this means that $Y_f$ is smooth projective and $D_f$ is a reduced anticanonical snc divisor.}, we call $(Y_f, D_f)$ a \emph{log Calabi--Yau compactification} of $f$. Removing the fiber over $\infty$ recovers the Calabi--Yau compactification $U_f$. 
\end{definition}
We make a few remarks on this definition:
\begin{itemize}
	\item The series $\pi_f(t)$ encodes period integrals along the fibers $f=\tfrac{1}{t}$, where $f$ is viewed as a regular function on $(\CC^\times)^n$. The period condition can then be interpreted as an equality of Gromov--Witten invariants on the ``A-side" and period integrals on the ``B-side", akin to formulations of mirror symmetry for Calabi--Yau threefolds.
	\item If $L=0$ is the trivial divisor class, the series $\tilde{I}^{(X, L)}(t)$ is also known as the regularized quantum period, denoted $\widehat{G}_X(t)$ in \cite{MirrorSymmetryFano}.
\end{itemize}
Extensive computational evidence suggests that toric Landau--Ginzburg models for $X=(X, 0)$ have very special coefficients. If $\dim(X) \leq 3$ and $-K_X$ is very ample, then \cite[Corollary 2]{Przyjalkowski} shows that $X$ has a toric Landau--Ginzburg model that is a \emph{Minkowski polynomial}. In view of the drawbacks of Minkowski polynomials lined out in~\cite[Section 4]{Coates}, the Fanosearch program~\cite{MirrorSymmetryFano} conjectures that toric Landau--Ginzburg models for Fano varieties are \emph{rigid maximally mutable Laurent polynomials} (rigid MMLPs). In particular, it is known that:
\begin{itemize}
	\item there is a one-to-one correspondence between smooth Fano threefolds with very ample anticanonical bundle and rigid MMLPs in three variables with reflexive Newton polytope up to an equivalence relation called \emph{mutation}~\cite[Theorem 4.1]{Coates}.
	\item every known Laurent polynomial satisfying the period condition for a Fano manifold $X$ in any dimension is a rigid MMLP (see~\cite[Section 5]{Coates} and references therein).
\end{itemize}
Suppose that $f$ is a Laurent polynomial in three variables with reflexive Newton polytope. 
To construct a Calabi--Yau compactification of a Laurent polynomial $f$, one can proceed as follows.\\
Let $Y_P$ be the Gorenstein toric Fano variety associated to the normal fan of $P$ (so that $P$ is the moment polytope of $Y_P$). Since $P$ is a reflexive 3D polytope, $Y_P$ admits a crepant resolution $\pi \colon (\bar{Y}, \bar{D}) \rightarrow (Y_P, D_P)$. Lattice points of $P$ give a basis of sections for $-K_{\bar{Y}}$. Therefore, $f$ and the constant polynomial $1$ define an anticanonical pencil on $\bar{Y}$, with induced rational map $\bar{Y}\dashrightarrow \PP^1$. The zero locus of the section $1$ is the toric boundary $\bar{D}$ of $\bar{Y}$, so that the base locus of the pencil is a union of (possibly nonreduced) curves on $\bar{D}$. The rational map can be resolved by an iterated blowup 
\[
Y_f \rightarrow \dots \rightarrow \bar{Y}
\]
Denoting $\tilde{f} \colon Y_f \rightarrow \PP^1$ the induced map, we have $\tilde{f}^{-1}(\infty)=D_f$.
In general, we cannot expect the iterated blowup $Y_f \rightarrow \bar{Y}$ to be a toric model in the sense of Definition~\ref{def:toricmodel}, as the curves in the base locus might be singular, or fail to have simple normal crossings with $1$-dimensional strata of $\bar{Y}$.
To remedy this, we must impose further conditions on the Laurent polynomial. Following \cite{Doran, CortiHackingPetracci}, we define:
\begin{definition}
	Let $P$ be a 3D reflexive lattice polytope. 
	\emph{Minkowski data} $M$ for $P$ assigns to each facet $F \subset P$ a decomposition
	\[
	F=k_1F_1+\dots +k_rF_r, \qquad k_i \in \ZZ_{> 0}
	\]
	into lattice polygons.
	An \emph{$M$-polynomial} is a Laurent polynomial $f$ with $\Newt(f)=P$ whose restriction to each facet $F$ factors as
	\[
	f|_F=f_1^{k_1}\dots f_r^{k_r}
	\]
	where the $f_i$ are Laurent polynomials with $\Newt(f_i)=F_i$. 
	We denote $\mathcal{L}_{P, M}$ the set of $M$-polynomials. \\
	We say that Minkowski data $M$ is \emph{admissible} if each $F_i$ is affinely equivalent to an $A_n$-triangle, that is, a lattice triangle with vertices $(0,0)$, $(n,0)$ and $(0,1)$ where $n \geq 0$. Note that an $A_0$-triangle is a line segment\footnote{\cite{CortiHackingPetracci} shift the index $n$ by $-1$.}. \\
	A \emph{Minkowski polynomial} is an $M$-polynomial for admissible $M$ such that after a change of variables,
	\begin{equation}\label{eq55}
		f_i=(1+x)^n+y.
	\end{equation}
\end{definition}
The base locus of the pencil defined by a $M$-polynomial can be easily described:
\begin{lemma}\label{lem:Minkowski}
	Let $f$ be a $M$-polynomial. Let $F \subset P$ be a facet with Minkowski decomposition 
	\[
	F=k_1F_1+\dots +k_rF_r
	\]
	For each $1 \leq i \leq r$, the base locus $\mathcal{B}$ on $D_F$ contains one curve $C_i$ of multiplicity $k_i$, of arithmetic genus equal to the number of interior lattice points of $F_i$. Moreover, $C_i$ does not contain any torus-invariant strata. 
	In particular, if $M$ is admissible, every curve in $\mathcal{B}$ is a smooth rational curve (possibly with multiplicity).
\end{lemma}
\begin{proof}
	Let $Z_f \subset {Y}_P$ be the anticanonical hypersurface defined by $f$. The part of the base locus $\mathcal{B}$ supported on the toric divisor $D_F$ is given by the vanishing of $f|_F$. By definition, we have a factorization $f|_F=f_1^{k_1} \dots f_r^{k_r}$, so that the base locus contains one curve $C_i:\{f_i=0\}$ of multiplicity $k_i$ for each $1 \leq i \leq r$. Standard results on toric surfaces (e.g. \cite[Theorem 10.5.8]{CoxLittleSchenck}) show that the arithmetic genus of $C_i$ is equal to the number of interior lattice points of $F_i$. Since $\Newt(f)=P$, the base locus does not contain any torus-invariant strata.\\
	If $M$ is admissible, each $F_i$ is an $A_n$-triangle, which is Minkowski indecomposable, so that $C_i$ is irreducible.  Furthermore, an $A_n$-triangle has no interior lattice points so that $p_a(C_i)=0$. Combined, this shows that $C_i$ is smooth rational.
\end{proof}
\begin{proposition}\cite[Proposition 3.6]{Doran}\label{prop:compactification}
	Let $P$ be a reflexive 3D lattice polytope. Let $(\bar{Y}, \bar{D}) \rightarrow (Y_P, D_P)$ be a crepant (toric) resolution. Let $M$ be Minkowski data for $P$ and let $f$ be a $M$-polynomial.
	There is an iterated blowup 
	\[
	(Y_f, D_f) \rightarrow \dots \rightarrow (\bar{Y}, \bar{D})
	\] 
	such that $f$ extends to a projective morphism $\tilde{f} \colon (Y_f, D_f) \rightarrow (\PP^1, \infty)$ with $U_f=Y_f \setminus D_f$ a Calabi--Yau compactification of $f$.
	Moreover
	\begin{itemize}
		\item If $M$ is admissible, then every blowup center is a smooth rational curve.
		\item if $f$ is very general, then every blowup center has simple normal crossings with the boundary and hence $D_f$ is snc. In particular, $(Y_f, D_f) \rightarrow (\bar{Y}, \bar{D})$ is a toric model and $(Y_f, D_f)$ is a log CY compactification of $f$.
	\end{itemize}
\end{proposition}
\begin{proof}
	Let $\mathcal{B}$ be the base locus of the anticanonical pencil on $Y_P$ defined by $f$. 
	Let $\bar{Y} \rightarrow Y_P$ be a crepant resolution.
	If $Y_P$ is singular along a $1$-dimensional stratum $D_e$, we denote $E^e_j$ the exceptional divisors dominating $D_e$, and $\Gamma^e_j$ a general fiber of $E_j^e \rightarrow D_e$, a smooth rational curve.
	Since $\mathcal{B}$ avoids zero-dimensional strata of $Y_P$ by Lemma~\ref{lem:Minkowski}, the strict transform $\mathcal{\bar{B}}$ consists of the curves in $\mathcal{B}$, possibly together with some of the fibers $\Gamma^e_j$, with some multiplicity. 
	Let $C_0$ be a component of $\mathcal{\bar{B}}$ of multiplicity $k$. Denote $\bar{D}_0$ the component of $\bar{D}$ containing $C_0$. Let $\tilde{Y} \rightarrow \bar{Y}$ be the blowup of $C_0^{red}$. By Lemma~\ref{lem:Minkowski} again, $C_0^{red}$ does not contain a torus-invariant stratum of $\bar{D}$, so the exceptional divisor $E$ has discrepancy $0$ and the strict transform $\tilde{D}$ of $\bar{D}$ is an anticanonical divisor. 
	Let $\tilde{\mathcal{B}}$ be the base locus of the strict transform of the pencil. We have that $E \cap \tilde{D}_0=C_0$, whereas $E \cap \tilde{D}_j$ for components $D_j$ adjacent to $D_0$ consists of smooth rational curves transverse to $D_e$ (see Figure~\ref{fig2} for an illustration). 
	Each of the curves in $E \cap \tilde{D}$ appears in $\tilde{\mathcal{B}}$ with multiplicity less than $k$; in particular, $\tilde{\mathcal{B}}=\mathcal{\bar{B}} \setminus C_0$ if $k=1$. \\
	We now successively blow up all curves in the base locus, starting with curves of highest multiplicity. Each round of blowup introduces finitely many new curves to the base locus, but reduces the maximal multiplicity. This ensures that the base locus can be resolved by a finite number of blowups $Y_f \rightarrow \dots \rightarrow \bar{Y}$. By construction, $f$ extends to $\tilde{f} \colon Y_f \rightarrow \PP^1$ with $D_f=\tilde{f}^{-1}(\infty)$ the strict transform of the toric boundary, an anticanonical divisor on $Y_f$.\\
	If $M$ is admissible, then every component of $\mathcal{{B}}$ is a smooth rational curve by Lemma~\ref{lem:Minkowski}. We argued above that each of the curves $\Gamma^e_j$ and $E \cap \tilde{D}_j$ is smooth rational as well, so that the base locus can be resolved by only blowing up along smooth rational curves. \\
	Let now $f$ be very general. Let $C_i$ be a component of $\mathcal{B}$, defined by a Laurent polynomial $f_i$. $C_i$ has simple normal crossings with the boundary if and only if $C_i$ is smooth, and the roots of $f_i$ along edges of $\Newt(f_i)$ are distinct. Since $f$ is very general, both conditions are clearly satisfied. Each of the curves $\Gamma^e_j$ and $E \cap \tilde{D}_j$ has simple normal crossings with the boundary as well, so that $(Y_f, D_f) \rightarrow (\bar{Y}, \bar{D})$ is a composite of interior blowups and hence a toric model, in the sense of Definition~\ref{def:toricmodel}. 
\end{proof}
\begin{remark}The construction of the log Calabi--Yau compactification $(Y_f, D_f)$ is not canonical, as it depends on the order in which the blowups are performed. However, any two log CY compactifications are related by a sequence of flops.
\end{remark}
We emphasize that the condition that the centers have simple normal crossings with the toric boundary often fails for Minkowski polynomials. The specialized form of the coefficients in \eqref{eq55} implies that if $F_i$ is an $A_n$-triangle and $n>1$, the curve $C_i:\{f_i=0\}$ is tangent to order $n$ to a one-dimensional stratum of $\bar{Y}$. In the blowup of $\bar{Y}$ in $C_i$, the strict transform of $\bar{D}_F$ has an $A_{n-1}$-singularity, so that $(Y_f, D_f)$ is not a log Calabi--Yau pair. However, in all known examples where such an $A_n$-triangle occurs, the base locus on the adjacent facet contains curves, all transverse to this stratum, which can be blown up first in the construction of the compactification. The strict transform of $C_i$ then becomes transverse to the stratum as well, and can subsequently be blown up. \\
We expect that it is possible to choose such an ordering of the curves for every Minkowski polynomial $f$ and therefore obtain a log CY compactification $(Y_f, D_f)$ as in Proposition~\ref{prop:compactification}, but are not aware of a general proof of this.
	
	\begin{proposition}\label{prop:mirrorfamily}
		Let $P$ be a reflexive 3D lattice polytope, let $M$ be admissible Minkowski data for $P$, and let $f$ be a very general $M$-polynomial.
	Let $\pi \colon (\cY, \cD) \rightarrow B$ be a family of log Calabi--Yau pairs with central fiber $(Y_f, D_f)$. 
	Then Conjecture~\ref{conj:torelli1} holds for very general fibers of the family: if $b_1$ and $b_2$ are very general closed points in $B$, any parallel transport isomorphism
	\[
	H^3(\mathcal{U}_{b_1}) \cong H^3(\mathcal{U}_{b_2})
	\]
	of MHS is induced by an isomorphism $(\cY_{b_1}, \cD_{b_1}) \cong (\cY_{b_2}, \cD_{b_2})$.
\end{proposition}
\begin{proof}
	We have to check that the assumptions of Theorem~\ref{thm:GenericTorelli} are satisfied. This follows immediately from Proposition~\ref{prop:compactification}: $(Y_f, D_f)$ admits a toric model since $f$ is very general, and since $M$ is admissible, every center in the toric model is a smooth rational curve. 
\end{proof}
Families of log Calabi--Yau pairs $\pi \colon (\cY, \cD) \rightarrow B$ with central fiber $(Y_f, D_f)$ as in Proposition~\ref{prop:mirrorfamily} are studied in \cite{Doran}: the authors show in \cite[Theorem 1.8]{Doran} that there is a family of Landau--Ginzburg models $(\cY, \cD, {\bf f})$ over a base $B$ of dimension $\rk\Pic(X)$. As $t \in B$ varies, general fibers of ${\bf f}_t$ form a family of $\Pic(X)^\vee$-polarized $K3$-surfaces with dominant period map. Here $\Pic(X)^\vee$ denotes the Dolgachev--Nikulin dual. \\
Under mirror symmetry, very general admissible $M$-polynomials as in Proposition~\ref{prop:mirrorfamily} are expected to be toric Landau--Ginzburg models for pairs $(X, L)$, where $X$ is a smooth Fano threefold with very ample anticanonical bundle, and $L \in \Pic(X) \otimes \CC$ is very general (see \cite[Section 1.4]{Doran}).

\bibliography{torelli3d}
\bibliographystyle{alpha}
\end{document}